\documentclass[16pt]{article}
\usepackage{graphicx}
\usepackage{amssymb}

\newtheorem{theorem}{Theorem}
\newtheorem{proposition}{Proposition}

\newtheorem{lemma}{Lemma}

\begin{document}

\title{Link projections and flypes}
\author{Cam Van Quach Hongler and Claude Weber}

\date{}


\maketitle
\begin{center}
  En souvenir de Felice Ronga
\end{center}
\begin{abstract}
Let $\Pi$ be a link projection in $S^2$.  John Conway \cite{Conway} and later Francis Bonahon and Larry Siebenmann \cite{BoSi} undertook to split $\Pi$ into canonical pieces. These pieces received different names: basic or polyhedral diagrams on one hand, rational, algebraic, bretzel, arborescent diagrams on the other hand. This paper proposes a thorough presentation of the theory, known to happy fews. We apply the existence and uniqueness theorem for the canonical decomposition to the classification of Haseman circles and to the localisation of the flypes.
\end{abstract}

\section{Introduction}

\subsection{The canonical decomposition of a link projection}

Let $\Pi$ be a link projection in $S^2$.  John Conway \cite{Conway} and later Francis Bonahon and Larry Siebenmann \cite{BoSi} undertook to split $\Pi$ into canonical pieces. These pieces received different names: basic or polyhedral diagrams on one hand, rational, algebraic, bretzel, arborescent diagrams on the other hand. 

The clue to break a projection $\Pi$ into pieces is to consider circles (simple closed curves) $\gamma$ embedded in $S^2$, cutting $\Pi$ transversally far from  double points. Let $\sharp(\gamma \cap \Pi)$ be the (even) number of intersection points between $\gamma$ and $\Pi$. We wish that circles $\gamma$ with $\sharp(\gamma \cap \Pi) \leq 2$ cut $\Pi$ in a trivial way. This is equivalent to require that $\Pi$ is connected (as a subset of $S^2$) and that $\Pi$ is indecomposable with respect to connected sum of projections (i.e. prime). 

{\bf Rule 1.}  The projections we consider are connected and prime. 

Circles $\gamma$ with $\sharp(\gamma \cap \Pi) = 4$ have been first considered by Mary Haseman in \cite{Haseman}. Pictures drawn in \cite{Tait} suggest that Peter Tait was close to this concept. John Conway made an extensive use of them, but we shall call them {\bf Haseman circles} for historical reasons. 

{\bf Definition.} A Haseman circle $\gamma$ is said to be {\bf compressible} if it bounds a disc $\Delta$ in $S^2$ such that there exists a properly embedded arc $\alpha \subset \Delta$ with $\alpha \cap \Pi = \emptyset$ and with $\alpha$  not boundary parallel. The arc $\alpha$ is called a {\bf compressing arc} for $\gamma$ in $\Delta$.

{\bf Comment.} Since we assume $\Pi$ to be prime, the only possible picture for a compressible $\gamma$ is as represented  in Figure 1.

\begin{figure}[ht]    
   \centering
    \includegraphics[scale=0.2]{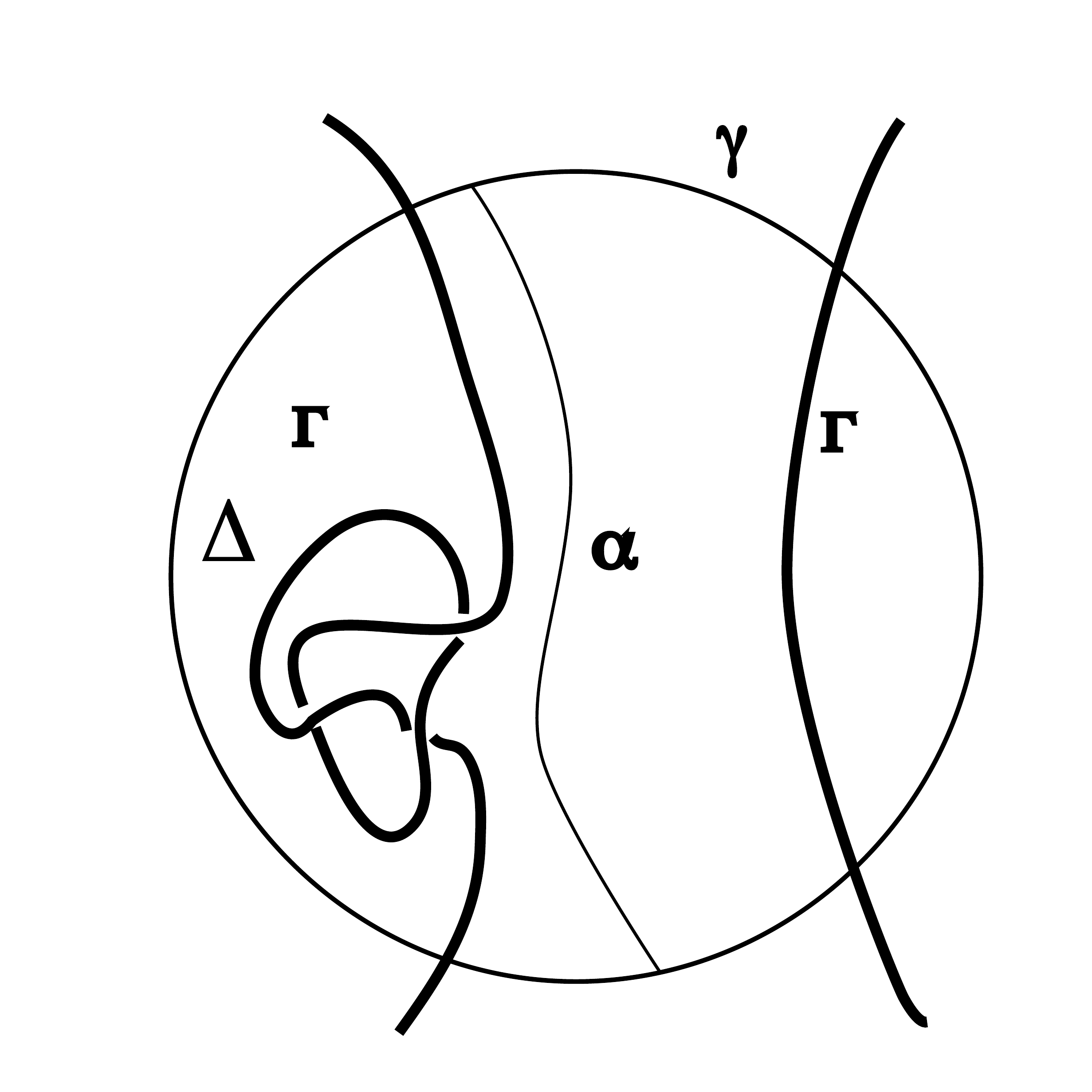}
\caption{ A  compressible Haseman circle}
\end{figure}

{\bf Rule 2.} Haseman circles are assumed  incompressible.

{\bf Definitions.} 1) Two Haseman circles $\gamma_0$ and $\gamma_1$ are {\bf isotopic} if there is a path $\gamma_t$ with $0 \leq t \leq 1$ joining them among Haseman circles.
\\
2) Two Haseman circles are {\bf parallel} if they bound an annulus $A$ such that the pair $(A , A \cap \Pi )$ is diffeomorphic to Figure 2.
\\
Disjoint Haseman circles are isotopic iff they are parallel ( Figure 2) .
\begin{figure}[ht]    
   \centering
    \includegraphics[scale=0.4]{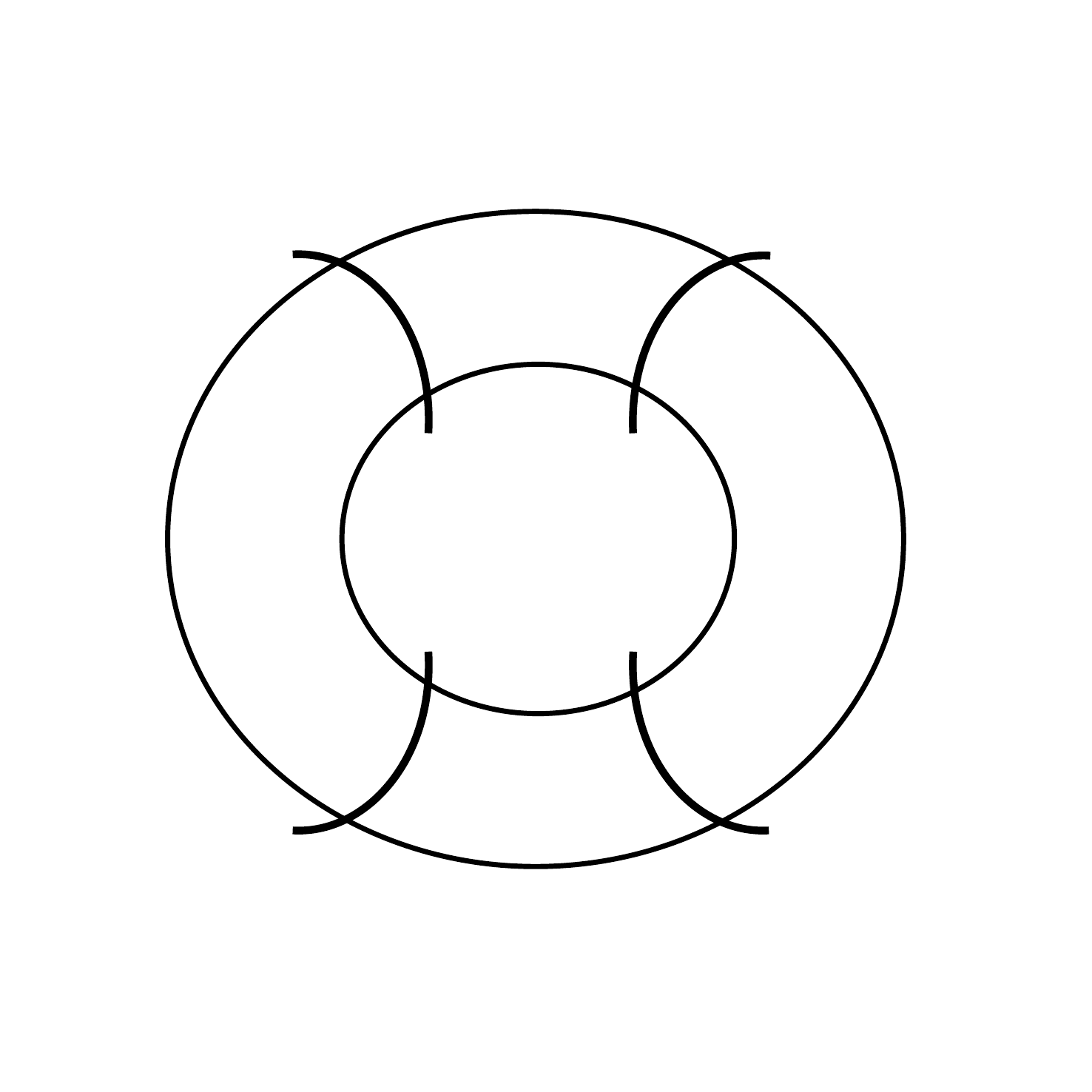}
\caption{ Isotopic Haseman circles}
\end{figure}

{\bf Definition.}  A {\bf family of Haseman circles} $\mathcal{H}$ for $\Pi$  is a collection of Haseman circles satisfying the following conditions.
\\
1. Any two circles are disjoint.
\\
2. No two circles are parallel.

Note that a family is always finite because a link projection has a finite number of crossing points.

Let $\mathcal{H} = \{ \gamma_1 , ... , \gamma_n \}$ be a family of Haseman circles for $\Pi$. Let $R$ be the closure of a connected component of $S^2 \setminus \bigcup_{i=1}^{i=n} \gamma_i$. We call the pair $(R , R \cap \Pi)$ a {\bf diagram}  determined by the family $\mathcal{H}$. We say more about diagrams in Section 2. Among diagrams, we shall single out:
\\
1) singletons;
\\
2) basic diagrams which contain few Haseman circles up to isotopy (Haseman circles are boundary parallel); 
\\
3) band diagrams which contain many. 

{\bf Definitions.} A family of Haseman circles is called a {\bf Conway family} if each diagram determined by this family is either a singleton or a basic diagram or a band diagram.  A Conway family is {\bf minimal} if the removal of any Haseman circle from this family transforms it into a non-Conway family. 

Broadly speaking, the {\bf Main Theorem} says that for a given projection $\Pi$, minimal Conway families exist and are unique up to isotopy. See Section 2 for a precise statement. The diagrams determined by this minimal family $\mathcal{C}_m$  are the elements of the canonical splitting mentioned in the beginning of this section. 

The Main Theorem has two important consequences which are presented in the next two subsections.

\subsection{Classification of Haseman circles} 

The key step in the proof of the Main Theorem (given in Section 3) is the following lemma.

{\bf Push-off lemma.} Let $\mathcal{C}_m$ be a minimal Conway family for the projection $\Pi$. Let $\gamma$ be some Haseman circle. Then there is an isotopy from $\gamma$ to $\gamma'$ such that $\gamma'  \cap \gamma_i = \emptyset$ for each member $\gamma_i$ of $\mathcal{C}_m$.  

The Haseman circle $\gamma'$ is contained in one of the diagrams determined by $\mathcal{C}_m$. Let $D$ be one of them.  In $D$, $\gamma$ may be isotopic to a boundary component $\beta$ of $D$, in other words to an element of $\mathcal{C}_m$. In this case there are in fact two diagrams into which $\gamma$ can be isotoped: $D$ and the diagram $D'$ which has also $\beta$ as a boundary component. If $\gamma$ is not isotopic to a member of $\mathcal{C}_m$ the diagram $D$ into which it can be isotoped is unique. Moreover in a basic or in band diagram, Haseman circles can be easily classified up to isotopy. As a consequence, from the knowledge of $\mathcal{C}_m$ we can understand all Haseman circles for $\Pi$. See also the end of Section 3.

\subsection{Flypes}

We recall the definition of a {\bf flype}, originally due to Tait under the name of {\bf distorsion}. Consider Figure 3a, representing a projection in $S^2$ with two Haseman circles.

The two discs $A$ and $B$ contain most of the diagram. They are bounded by the two Haseman circles. It is easy to see that the projection can be transformed  by a 3-dimensional isotopy to the  one represented in Figure 3b.

Disc $B$ is unchanged. Disc $A$ is altered by a rotation of angle $\pi$ and axis contained in the plane of the figure. The disc $A$ (and its content) is called the {\bf moving disc}. The crossing point $P$ is called {\bf active}. Clearly, there is a similar transformation for which $B$ is moving. This  transformation was discovered and intensively used by Peter Tait in \cite{Tait}. Since the publication of John Conway's famous paper \cite{Conway} it is now called a flype.

\begin{figure}[ht]    
   \centering
    \includegraphics[scale=0.45]{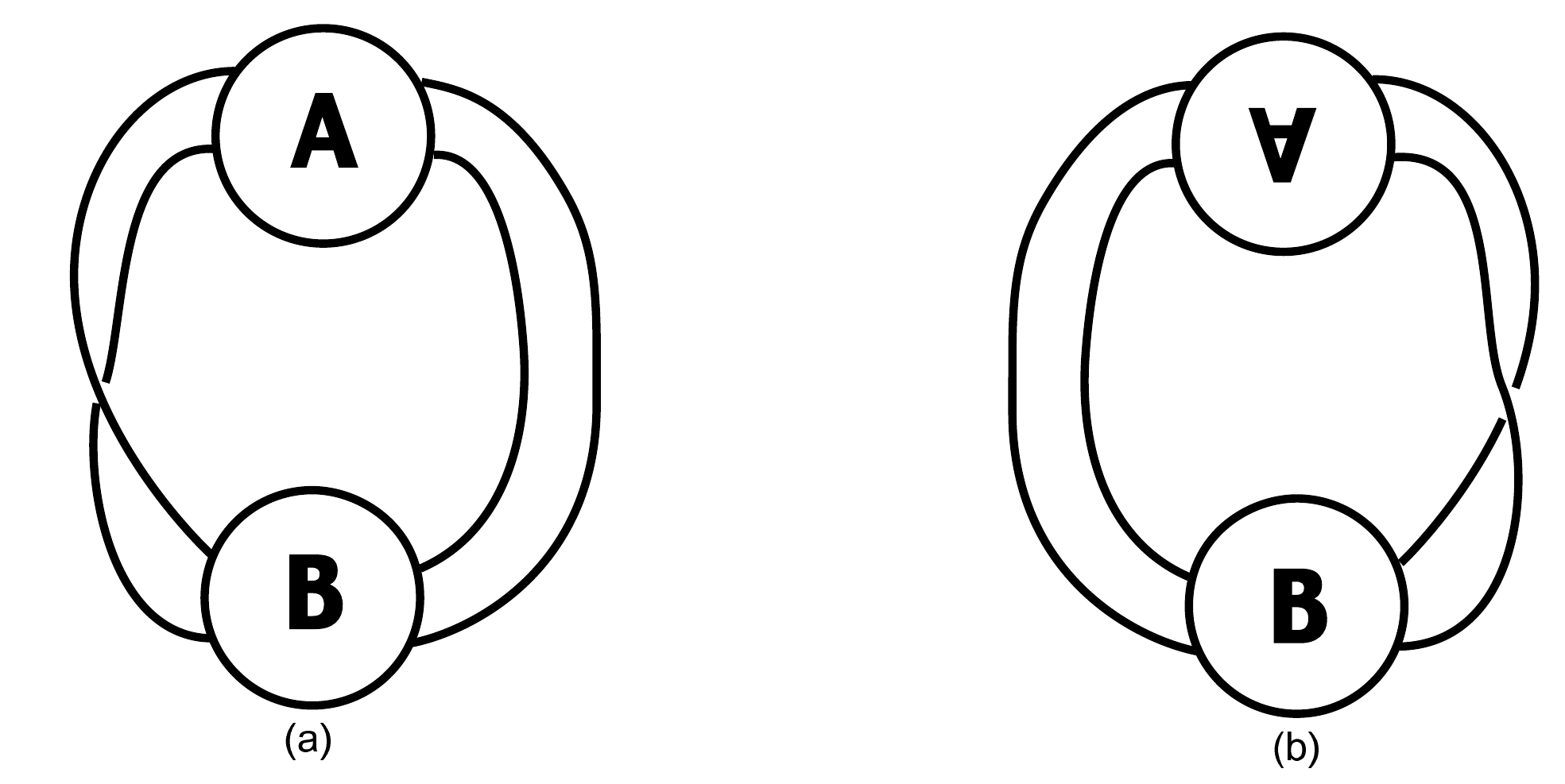}
\caption{Before and after the flype}
\end{figure}

{\bf Remark.}  As we assume that the boundary of disc $A$ and of disc $B$ are Haseman circles, they are incompressible and hence there are crossing points inside both discs. Thus, our point of view is slightly different from that of Menasco-Thistlethwaite in \cite{MeTh}.

William Menasco and Morwen Thistlethwaite prove in \cite{MeTh} the {\bf Flyping Theorem}:  any two minimal alternating prime diagrams which represent the same oriented link type differ by a finite sequence of flypes. This theorem has important consequences for the classification of alternating knots and links and for the study of their symmetries. We plan to devote a future publication  to this second topic. In this context, the following question is natural:

{\bf Question.} How can we find the position of flypes?

Specialists know that the answer is contained in the book by Francis Bonahon and Larry Siebenmann  (alas not yet published) \cite{BoSi}. The answer to the question is  the following:

1) flypes are localized in the arborescent part of the diagram;

2) the planar trees (one tree for each connected component of the arborescent part) which code the arborescent part give an excellent account of the position of flypes.

We offer here an elementary and self-contained  presentation of what is needed from Bonahon-Siebenmann to capture the flypes. 

It is interesting to compare our paper with \cite{Calvo} where the author takes a different point of view.

Here is an application. Consider an alternating projection
$\Pi$, satisfying Rule 1 representing a knot/link $K$ in $S^3$. A
consequence of Menasco and Thistlethwaite Flyping Theorem is that the
elements of the group $Sym(K)$ of symmetries of $K$ are compositions of
flypes and symmetries of the projection $\Pi$. In principle $Sym(K)$ can
hence be determined  ``by hand" (i.e. without a computer) and the
operation is made simpler if we know where the flypes are. If there are
no flypes for $\Pi$ then $Sym(K)$ is isomorphic to the group $Sym(\Pi)$
of symmetries of the projection.  These are sometimes called ``visible
symmetries". In other words, when there are no flypes, the visible
symmetries are the only symmetries. Compare with the book of Akio Kawauchi
\cite{Kawauchi} chapter 10, written by Makoto Sakuma. An easy consequence
of the arguments presented in Sections 4 and 5 is that there are no
flypes if and only if for each planar tree representing an arborescent
component of $\Pi$ the vertices of valency $\geq 2$ have weight 0.

Here are two examples built from Conway's polyhedron 6*.  Now 6* is made
of the edges of an octahedron embedded in $S^2$. See Figure 9. Hence the group of
symmetries of 6* is isomorphic to the group $S_4 \times C_2$ of isometries
of the octahedron. If we replace the singletons of 6* by more complicated
diagrams, the symmetry group usually diminishes heavily.

Consider Figures 4 and 5. They are obtained from 6* by replacing as indicated some singletons by arborescent diagrams. Both projections are alternating without flypes (for projections without complications as these ones there is no need of a big machinery to see that. Our examples are presented only to illustrate the general principles).  Figure 4 represents a knot which is only -achiral (symmetry group $D_1$, just one reflection). Figure 5 represents a 3-component link with trivial symmetry group.

\begin{figure}[ht]    
   \centering
    \includegraphics[scale=0.4]{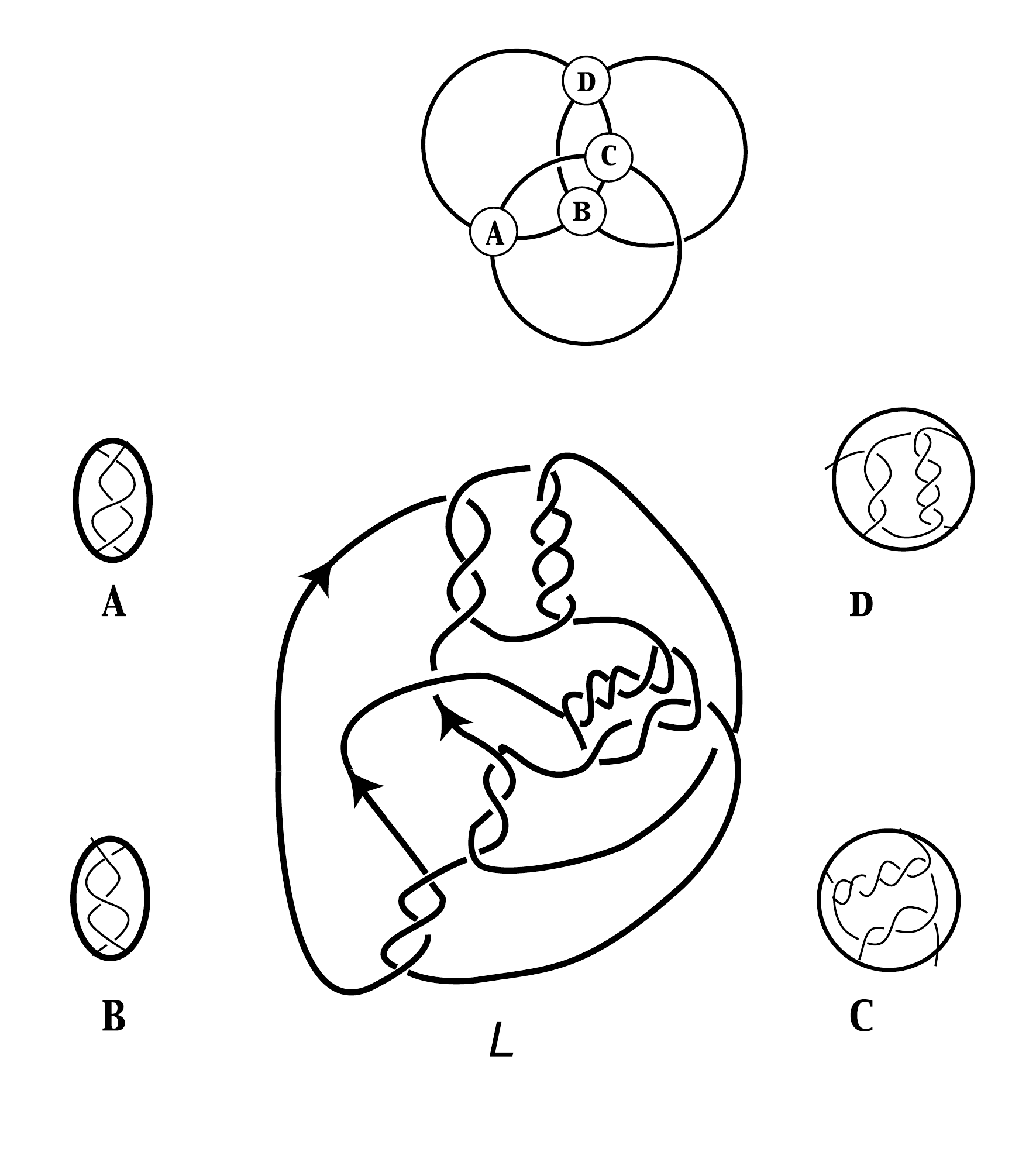}
\caption{A negatively achiral knot}
\end{figure}

\begin{figure}[ht]    
   \centering
    \includegraphics[scale=0.4]{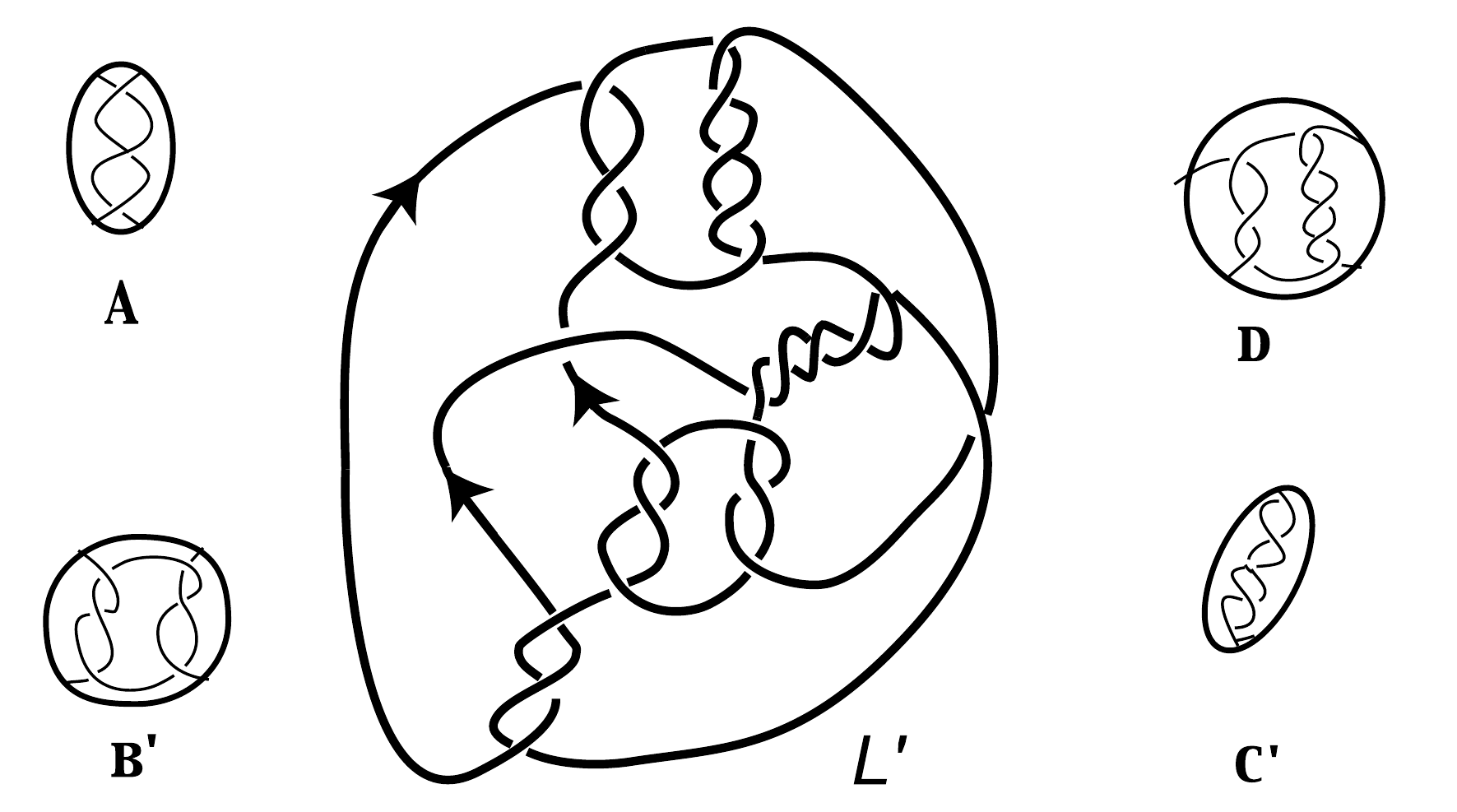}
\caption{ A link with trivial symmetry group}
\end{figure}

\noindent

\subsection{Content of the paper}

The content of our paper is as follows. In  Section 2, we define  diagrams and the three kinds of diagrams we are particularly  interested in: singletons, basic and band diagrams. We  consider  Conway families of Haseman circles. We state the Main Theorem (existence and uniqueness  of minimal Conway families). We define  the arborescent and the polyhedral part of a link projection. In Section 3, we prove the main theorem. A key ingredient is the push-off lemma. In Section 4, we use the main theorem and the push-off lemma to locate the flypes, as promised. In Section 5, we present Bonahon-Siebenmann's clever and extremely useful method for encoding an arborescent diagram by a planar tree (hence the name ``arborescent"). In Section 6, we show how a flype acts on the planar tree (again, this is known to Bonahon-Siebenmann). In Section 7, we apply the theory to alternating diagrams and links. In Section 8, a practical guide for the construction of the minimal Conway family and the planar tree is suggested.

A few words about our vocabulary are necessary. We have decided (somewhat reluctantly) not to use the word ``tangle", for several reasons. The main one is that Bonahon-Siebenmann use it to describe a 3-dimensional object. ``Tangle diagram" would have been acceptable, but then ``diagram" is shorter, and this is the term we shall use. Historically, it seems that the name tangle was first introduced by Mary Haseman in \cite{Haseman}. Therefore, we call ``Haseman circle" what is often called a ``Conway circle". A key contribution made by John Conway (``algebraic" vs ``polyhedral") is remembered here by the concept of ``Conway family of Haseman circles".

We emphasize the fact that our setting is always 2-dimensional.

{\bf Acknowledgements.} We thank Francis Bonahon and Larry Siebenmann for writing \cite{BoSi} and for making unpublished texts available to us . We really hope that their book will eventually be published~! Also special thanks to Michel Boileau. We thank Morwen Thistlethwaite for his much appreciated visit to Geneva and for having said to us one morning: ``look at the arborescent part". We also thank John Steinig for his precious help to improve our English.

We thank the Fonds National Suisse de la Recherche Scientifique for its support. The first author wishes to thank the organizers of the Quantum Topology Conference Thang Le and Stavros Garoufalidis   and  the NSF for its  financial support.

\section{Main Theorem}

\subsection{Diagrams}

A {\bf planar surface} $\Sigma$ is a compact connected surface embedded in the 2-sphere $S^2$. We denote by $v$ the number of connected components of the boundary $b\Sigma$ of $\Sigma$. We consider compact graphs $\Gamma$ embedded in $\Sigma$ and satisfying the following four conditions.

1) Vertices of $\Gamma$ have valency 1 or 4. 

2) Let $b\Gamma$ be the union of the vertices of $\Gamma$  of valency 1. Then $\Gamma$ is properly embedded in $\Sigma$, i.e. $b\Sigma \cap \Gamma = b\Gamma$.

3) The number of vertices of $\Gamma$ contained in each connected component of $b\Sigma$ is equal to 4.

4) A  vertex of $\Gamma$ of valency 4 is called a {\bf crossing point}.  We require that at each crossing point an over and an under thread be chosen and pictured as usual. We denote by $c$ the number of crossing points.
 
{\bf Definition.} The pair $D = (\Sigma , \Gamma)$ is called a {\bf diagram}.

{\bf Definition.} A {\bf singleton} is a diagram  diffeomorphic to Figure 6.

\begin{figure}[ht]    
   \centering
    \includegraphics[scale=0.5]{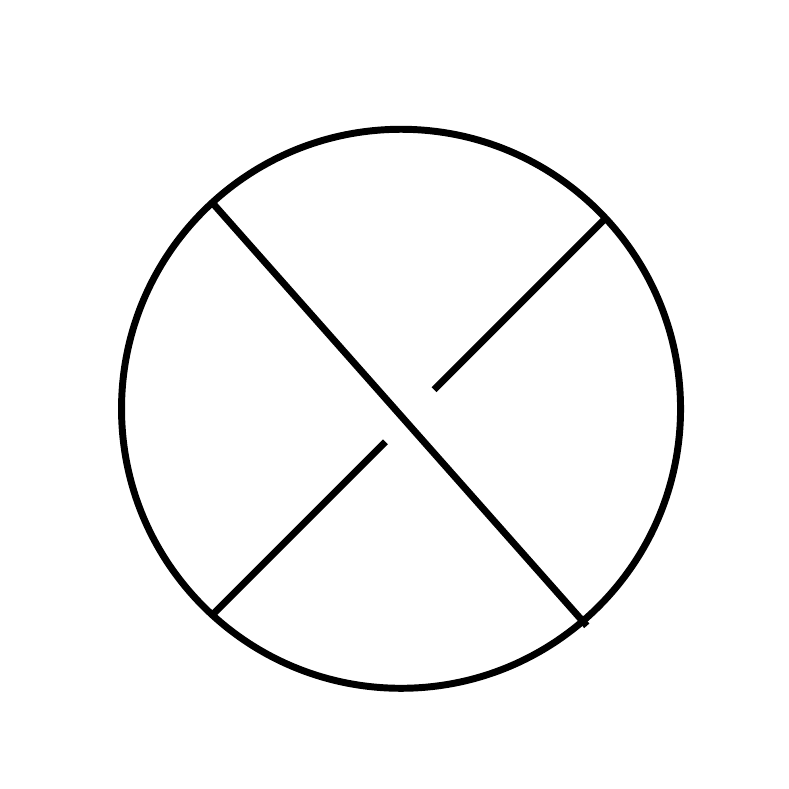}
\caption{ A singleton}
\end{figure}

 {\bf Definition.} A {\bf band diagram} is a diagram diffeomorphic to Figure 7 with $v \geq 3$. It can also be represented as in Figure 8.

\begin{figure}[ht]    
   \centering
    \includegraphics[scale=0.3]{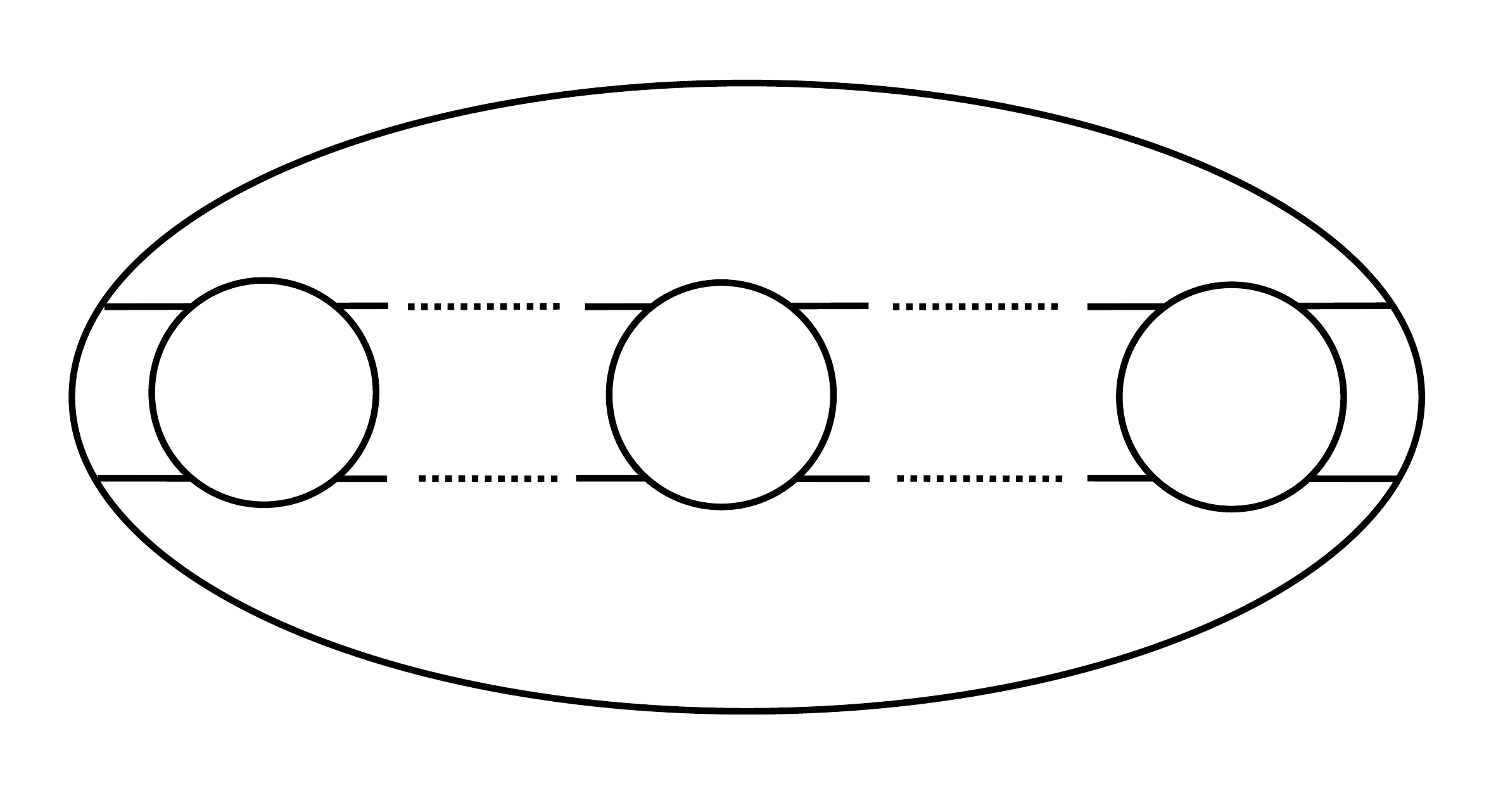}
\caption{ A band diagram}
\end{figure}

\begin{figure}[ht]    
   \centering
    \includegraphics[scale=0.3]{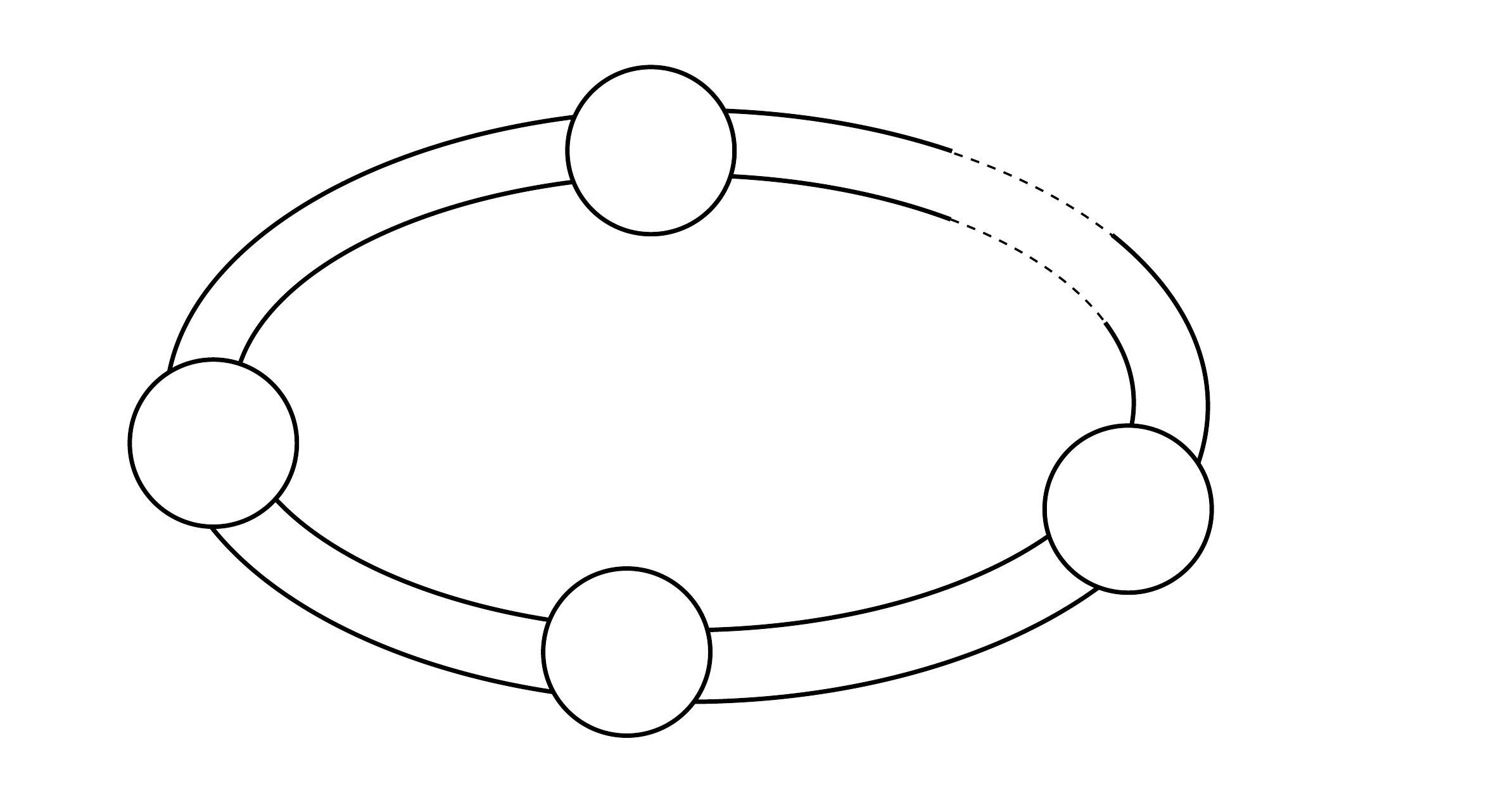}
\caption{A band diagram}
\end{figure}

{\bf Definition.} A diagram $D$ is {\bf basic} if:
\\
1) $D$ is not a singleton, nor a band diagram with $v = 3$, nor the diagram with $v = 1$ and $c = 0$.
\\
2) Every Haseman circle in $D$ is boundary parallel.

{\bf Comment.} Condition 2) is satisfied by the diagrams listed in 1) and we wish them to be excluded from basic diagrams. Basic diagrams are obtained from Conway polyhedra (indeed those which are indecomposable for tangle sum) by removing a small disc around each crossing point. Figure 9 is a picture of the simplest basic diagram, obtained from Conway $6^*$.

\begin{figure}[ht]    
   \centering
    \includegraphics[scale=0.3]{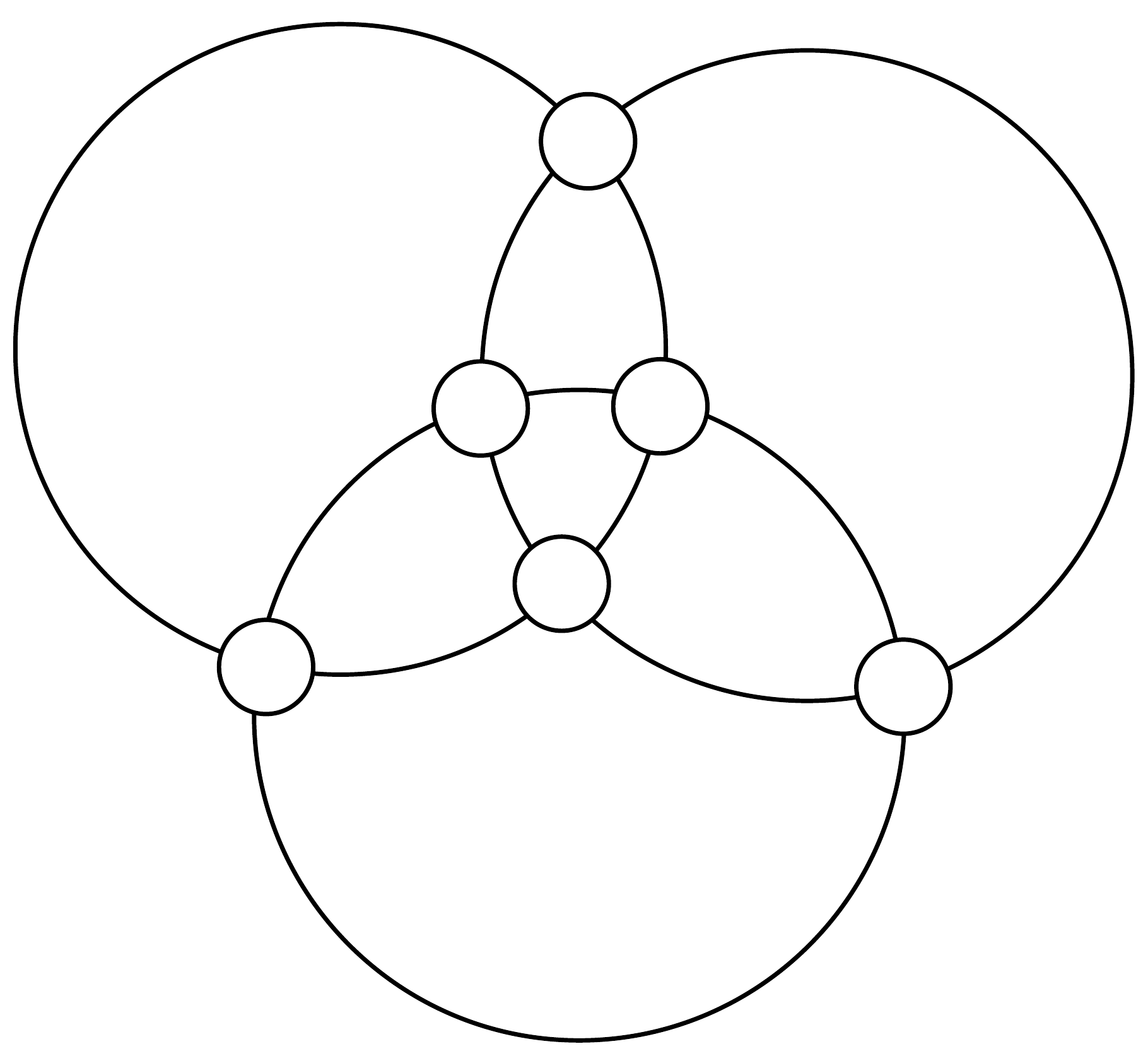}
\caption{ A basic diagram}
\end{figure}

\subsection{Haseman circles, Conway families and the Main Theorem}

{\bf Definition.} A {\bf Haseman circle} for a link projection $\Pi$ in $S^2$ is a circle $\gamma \subset S^2$ meeting $\Pi$ transversally in four points, far from  crossing points. A Haseman circle is said to be {\bf compressible} if
\\
i) $\gamma$ bounds a disc $\Delta$ in $\Sigma$.
\\
ii) There exists a properly embedded arc $\alpha \subset \Delta$ such that $\alpha \cap \Pi = \emptyset$ and such that $\alpha$ is not boundary parallel. The arc $\alpha$ is called a {\bf compressing arc} for $\gamma$.

{\bf Rule 2.} Haseman circles are  incompressible.

Two Haseman circles are said to be {\bf parallel} if they bound an annulus $A \subset \Sigma$ such that the pair $(A , A \cap \Gamma)$ is diffeomorphic to Figure 2.

Analogously, we define a Haseman circle $\gamma$ to be {\bf boundary parallel} if there exists an annulus $A \subset \Sigma$ such that:
\\
1) The boundary $bA$ of $A$ is the disjoint union of $\gamma$ and a boundary component of $\Sigma$;
\\
2) $(A , A \cap \Gamma)$ is diffeomorphic to Figure 2.

{\bf Definition.} Let $\Pi$ be a link projection. A {\bf family of Haseman circles} for $\Pi$ is a set of Haseman circles satisfying the following conditions.
\\
1. Any two circles are disjoint.
\\
2. No two circles are parallel.

Note that a family is always finite, since a projection has a finite number of crossing points.

Let $\mathcal{H} = \{ \gamma_1 , ... , \gamma_n \}$ be a family of Haseman circles for $\Pi$. Let $R$ be the closure of a connected component of $S^2 \setminus  \bigcup_{i=1}^{i=n} \gamma_i$. We call the pair $(R , R \cap \Gamma)$ a {\bf diagram} of $\Pi$ determined by the family $\mathcal{H}$.

{\bf Definition.} A family $\mathcal{C}$ of Haseman circles for $\Pi$ is called a {\bf Conway family} if each diagram determined by $\mathcal{C}$ is either:
\\
a) a singleton,
\\
b) a band diagram,
\\
c) a basic diagram.

Notice that neither a band diagram nor a basic diagram contains crossing points. Hence, if we have a Conway family, each crossing point is contained in a singleton.

{\bf Definition.} A Conway family $\mathcal{C}$ for $\Pi$ is {\bf minimal} if the deletion of any circle from $\mathcal{C}$ transforms it into a family which is no longer a Conway family.

The following is the MAIN THEOREM of the paper.

\begin{theorem}
(Existence and uniqueness of minimal Conway families.) Let $\Pi$ be a link projection in $S^2$. Then:
\\
i) There exist minimal Conway families for $\Pi$.
\\
ii) Any two minimal Conway families are isotopic, by an isotopy which respects $\Pi$.
\end{theorem}

The proof of the theorem will be given in Section 3.  The hypotheses we use are that the projection is prime and that Haseman circles are incompressible (Rules 1 and 2 above). $\Pi$ is not assumed to be alternating.

\subsection{Polyhedral and arborescent components}

{\bf Notation.} We write $\mathcal{C}_m$ for the {\bf minimal Conway family} of a link projection.

{\bf Definition.} A union of diagrams of $\Pi$ determined by  $\mathcal{C}_m$ is a {\bf pre-jewel} if it is:
\\
1) a connected union of one basic diagram and some singletons;
\\
2) maximal for these properties.

We now consider singletons. Suppose that we have a singleton whose boundary circle is also a boundary circle of a band diagram determined by $\mathcal{C}_m$. Using the local band, we can attribute a sign to the crossing point of the singleton by the rule pictured in Figure 10.

\begin{figure}[ht]    
   \centering
    \includegraphics[scale=0.6]{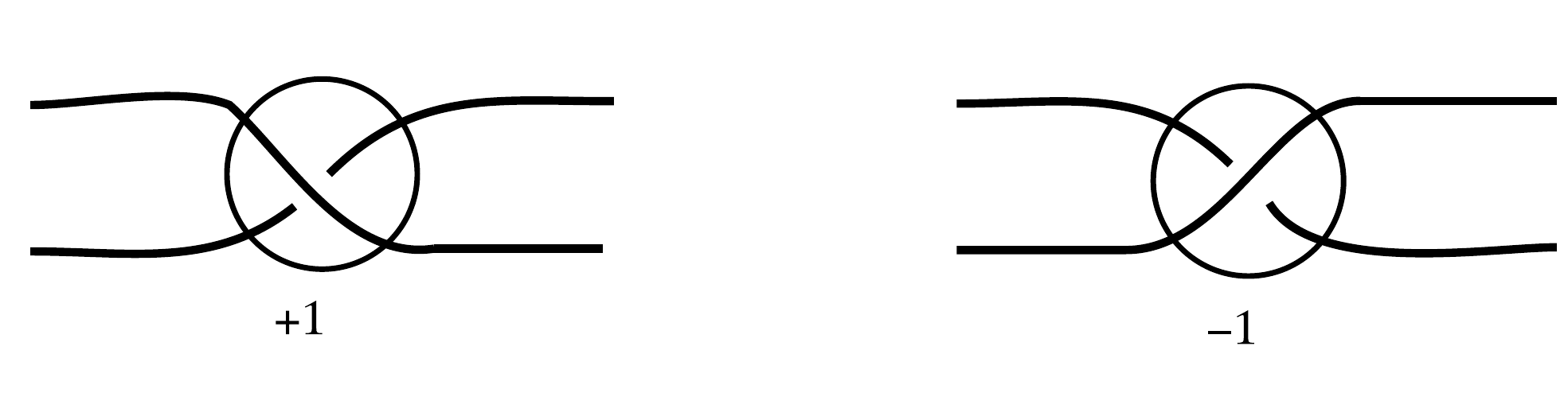}
\caption{Sign of the singleton}
\end{figure}

Note that singletons contained in a pre-jewel have no sign.

{\bf Rule 3.} Singletons sitting side by side along the same band have the same sign. In other words we assume that Reidemeister Move 2 cannot be applied to reduce the number of crossings.

Our next step is to get rid of singletons. For this, we remove from the family $\mathcal{C}_m$ the Haseman circles which bound a singleton. We obtain a new family $\mathcal{C}_{can}$ which we call the {\bf canonical Conway family}. We then examine the diagrams  determined by $\mathcal{C}_{can}$. They can be of two types: a former pre-jewel or a former band diagram with the Haseman circles bounding singletons removed. In the first case, we call the subdiagram a {\bf jewel}. In the second case, the subdiagram is called a {\bf twisted band diagram}. They are the canonical building blocks of $\Pi$. A typical twisted band diagram is pictured in Figure 11.

\begin{figure}[ht]    
   \centering
    \includegraphics[scale=0.4]{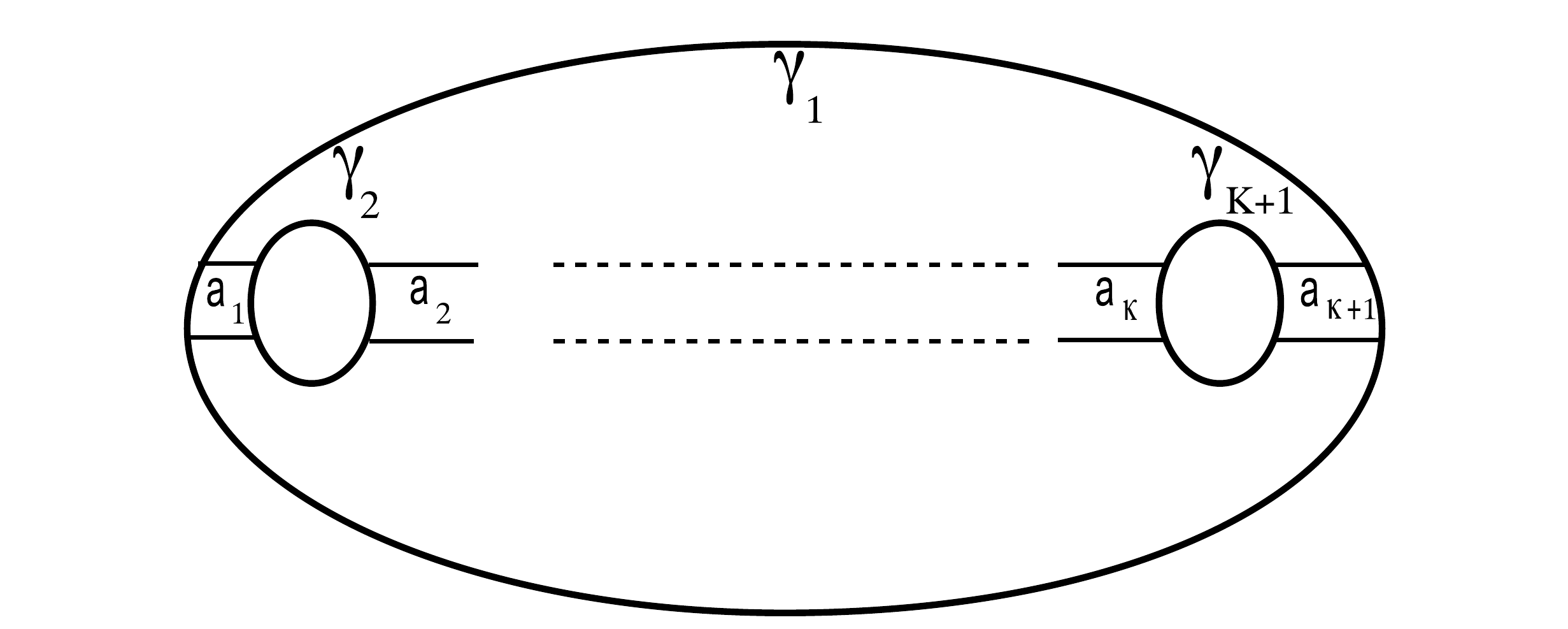}
\caption{A twisted band diagram}
\end{figure}

Here $\gamma_1 , ... , \gamma_{k+1}$ are canonical Haseman circles and we have $k \geq 0$. The $a_i$ are integers. $\mid a_i \mid$ denotes the number of former singletons sitting side by side between $\gamma_{i-1}$ and $\gamma _i$. The sign of $a_i$ tells us the sign of these singletons. The integer $a_i$ will be called an {\bf intermediate weight}. The corresponding portion of the diagram is called a {\bf twist}.

{\bf Properties of  intermediate weights.}
\\
i) It may happen that $k = 0$. In this case, we have $\mid a_1 \mid \geq 2$ because the diagram is obtained from a former band diagram.
\\
ii) If $k = 1$ then  $a_1 \neq 0$ or $a_2 \neq 0$.
\\
iii) If $k \geq 2$ it is possible that all $a_i$ are equal to zero.

{\bf Remark.} Using flypes and then Reidemeister Move 2, we can reduce the number of crossings of a twisted band diagram  in such a way that:
\\
1) either $a_i \geq 0$ for all $i = 1 , ... , k+1$;
\\
2) or $a_i \leq 0$ for all $i = 1 , ... , k+1$.

The reduction procedure is not quite canonical, but any two diagrams reduced in this manner are equivalent by flypes. This is enough for our purpose.

{\bf Rule 4.} We assume that in any twisted band diagram, all the $a_i$ which are different from zero have  the same sign.

We set  $a = \sum _{i=1}^{k+1} a_i$ and call $a$ the {\bf weight of the twisted band diagram}.

{\bf Remark.} It may occur that $\mathcal{C}_{can}$ is empty as in Figure 12.

\begin{figure}[ht]    
   \centering
    \includegraphics[scale=0.3]{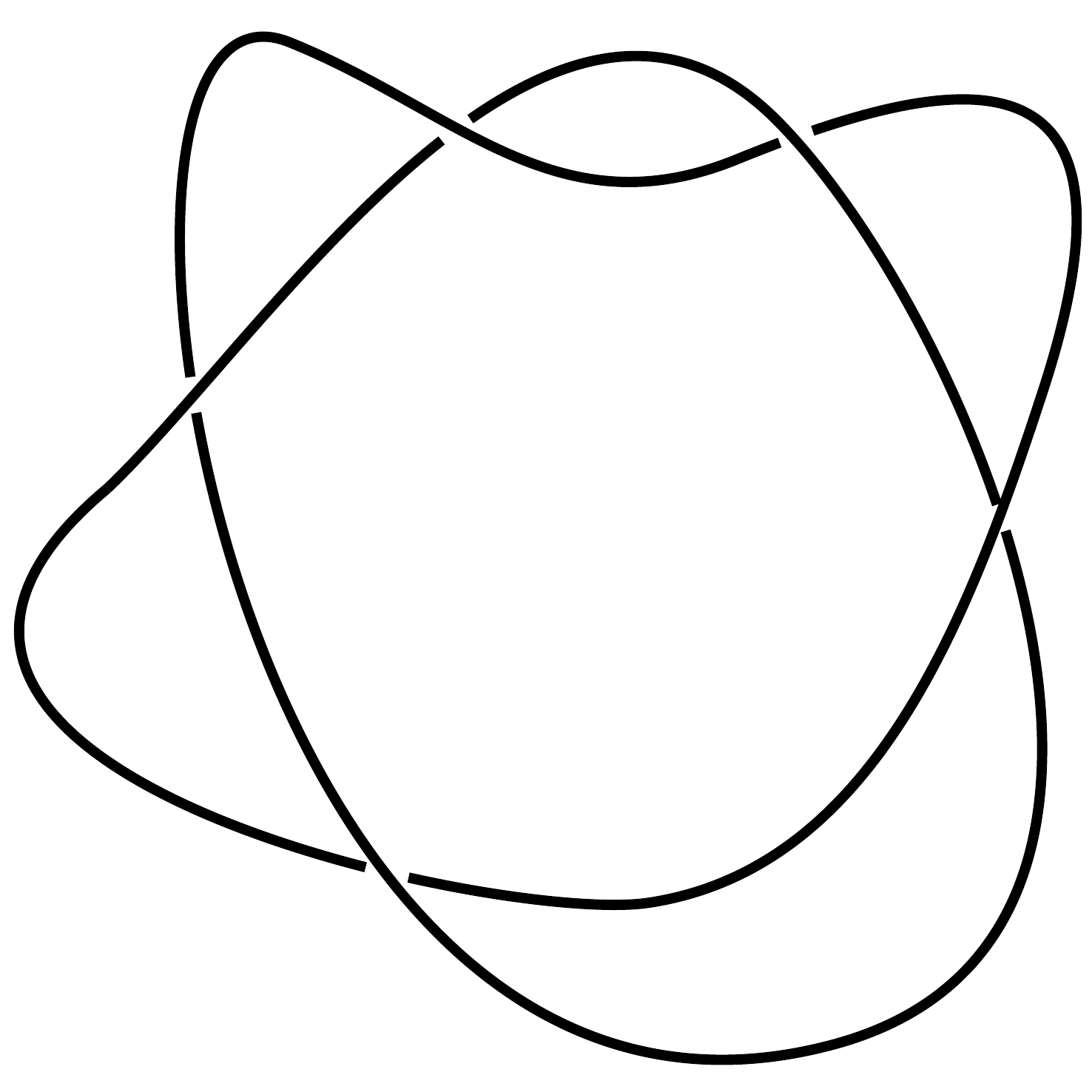}
\caption{A diagram without canonical circle}
\end{figure}

We are now ready to define polyhedral and arborescent components. Consider the canonical family $\mathcal{C}_{can }$ for the projection $\Pi$. Consider the diagrams determined by 
$\mathcal{C}_{can }$. These diagrams are either jewels or twisted band diagrams. 

{\bf Definitions.} A {\bf polyhedral component} of $\Pi$ is a connected union of jewels, maximal for inclusion. An {\bf arborescent component} of $\Pi$ is a connected union of twisted band diagrams, maximal for inclusion. The {\bf polyhedral part} is the disjoint union of the polyhedral components. The {\bf arborescent part} is the disjoint union of the arborescent components. One of these parts may be empty. In this case, the projection is said to be arborescent or polyhedral.

\section{Proof of the main theorem}

\subsection{Existence}

In this section, under- and over-passes at crossing points are unnecessary. Hence $\Gamma$ is just a graph embedded in $S^2$ satisfying conditions 1, 2, 3  from the beginning of Section 2. We require that Rules 1 and 2 be satisfied: graphs embedded in $S^2$ are prime and connected and Haseman circles are incompressible. 

We first prove that minimal Conway families exist. To do this, consider a maximal family $\mathcal{H}_M$ of Haseman circles. To be maximal means that we cannot add a new circle $\gamma$ (non-intersecting members of $\mathcal{H}_M$) to $\mathcal{H}_M$ in such a way that we still have a family. This means that $\gamma$ must be parallel to a member of $\mathcal{H}_M$. As $\Gamma$ has a finite number of crossing points, maximal families do exist. In general, they are not unique. We shall see why at the end of this section.

We claim that $\mathcal{H}_M$ is automatically a Conway family. First, let us remark that because $\mathcal{H}_M$ is maximal, each crossing point of $\Gamma$ is encircled by a member of  $\mathcal{H}_M$. Then, let us consider the diagrams of $\Gamma$ determined by  $\mathcal{H}_M$. Let $D'$ be one of them. From the definition of basic diagrams given in Section 2, we see that maximality implies that $D'$ is either a singleton, or a band diagram with three boundary components, or is basic. This concludes the proof of the existence of Conway families. From any such family, we can obtain a minimal one, by removing unnecessary circles.

\subsection{The push-off lemma}

The following lemma is crucial.

\begin{lemma}
(Push-off lemma.) Let $\Gamma$ be a graph in $S^2$ and let $\mathcal{C}_m$ be a minimal Conway family for $\Gamma$. Let $\gamma$ be some Haseman circle in $\Gamma$. Then we can isotope $\gamma$ to $\gamma '$, respecting $\Gamma$, in such a way that $\gamma '$ does not intersect any member of the family $\mathcal{C}_m$.
\end{lemma}

{\bf Proof of Lemma 1}

{\bf Step 1.} We move $\gamma$ by a small isotopy in such a way that it cuts  each circle of $\mathcal{C}_m$ transversally. We can also suppose that $\gamma$ intersects $\Gamma$ far from $\mathcal{C}_m$.

{\bf Step 2.} Each circle of $\mathcal{C}_m$ bounds two discs in $S^2$. Let $E$ be the set of  discs determined by each circle of $\mathcal{C}_m$ and let $\hat {E}$ be the subset of $E$ which consists of the discs which intersect $\gamma$. Let $\Delta$ be an element of $\hat {E}$ which is innermost among elements of $\hat {E}$. Write $C$ for the boundary of $\Delta$. ``$\Delta$ is innermost among elements of $\hat {E}$" means that $\gamma \cap C \neq \emptyset$ but that $\gamma \cap C' = \emptyset$ for any other element $C' \in \mathcal{C}_m$ contained in $\Delta$. 
\\
Let $\alpha$ be a connected component of $\gamma \cap \Delta$. It is a properly embedded arc in $\Delta$. Write $i_{\alpha}$ for the cardinality of $\alpha \cap \Gamma$. Since $\gamma$ cuts $\Gamma$ in 4 points, we have $0 \leq i_{\alpha} \leq 4$.

{\bf Step 3.} Suppose that  $i_{\alpha} = 0$. Since $C$ is incompressible, $\alpha$ can be isotoped to an arc in $C$ without passing through $\Gamma$. Hence $\alpha$ can be pushed off $\Delta$.

{\bf Step 4.} Suppose that $i_{\alpha} = 1$. The extremities of $\alpha$ separate $C$ into two arcs $C_1$ and $C_2$. Since card$(C \cap \Gamma) = 4$ and since $i_{\alpha} = 1$, we have card$(C_1 \cap \Gamma) = 1$ and card$(C_2 \cap \Gamma) = 3$ up to renumbering $C_1$ and $C_2$. Now $C_1 \cup \alpha$ is a circle in $S^2$ which cuts $\Gamma$ transversally in two points: one in $\alpha$ and one in $C_1$. Since $D$ is assumed to be prime, $\alpha$ can be isotoped to $C_1$ (respecting $\Gamma$) and hence off $\Delta$.

{\bf Step 5.} Suppose that $i_{\alpha} = 3$ or 4. The previous arguments show that $\gamma$ can be pushed recursively off any element of $E$ which is distinct from $\Delta$. Hence $\gamma$ can be pushed into the interior of $\Delta$. Since $\Delta$ is innermost in $\hat {E}$, we see that now $\gamma$ does not intersect any element of $\mathcal{C}_m$.

{\bf Step 6.} Suppose that $i_{\alpha} = 2$. Break  $C$ up again as $C = C_1 \cup C_2$ with $C_1 \cap C_2 = \alpha \cap C$. Because $i_{\alpha} = 2$ we have that card$(C_j \cap \Gamma)$ is even for $j = 1, 2$. But we cannot have card$(C_j \cap \Gamma) = 0$ for some $j$, otherwise $C_j$ would be a compressing arc for $\gamma$. Hence card$(C_j \cap \Gamma) = 2$ for $j = 1, 2$.

{\bf Step 7.} From the arguments already presented, we see that we are left with the following situation:
\\
i) There are exactly two innermost discs $\Delta _1$ and $\Delta _2$ in $\hat {E}$ such that card$(\gamma \cap \Delta _j \cap \Gamma) = 2$ for $j = 1, 2$.
\\
ii) $\gamma \cap \Delta _j$ consists of a single arc $\alpha _j$ and that arc is not boundary parallel in $\Delta _j$.

{\bf Step 8.} Let $\Delta$ be either $\Delta _1$ or $\Delta _2$. Consider the diagrams of $\Gamma$ determined by $\mathcal{C}_m$. Let $D'$ be the subdiagram which contains the arc $\alpha = \gamma \cap \Delta$. $D'$ exists because $\Delta$ is innermost. The boundary $C$ of $\Delta$ is one of the boundary components of $D'$.
\\
We claim that we may suppose that   $D'$ is neither a singleton nor a  basic diagram. A singleton is excluded because $\alpha$ would be boundary parallel. Assume therefore that  $D'$ is basic. Then $\alpha \cup C_j ~~~(j = 1, 2)$ is a circle which cuts $\Gamma$ in four points (two in $\alpha$ and two in $C_j$, see Step 6). If $\alpha \cup C_j$ is incompressible, then it is a Haseman circle. Because $D'$ is basic, this circle is boundary parallel. And therefore $\alpha$ can be pushed off $\Delta$. Since $C$ is incompressible, $\alpha \cup C_1$ or  $\alpha \cup C_2$ is incompressible.

{\bf Step 9.} Thus we can add to Conditions i) and ii) the condition:
\\
iii) The diagram of $\Gamma$ determined by $\mathcal{C}_m$ which contains $\alpha _j$ is a band diagram $(j = 1, 2)$.

{\bf Claim.} Such a situation is impossible since $\mathcal{C}_m$ is assumed to be minimal. 

Clearly, the claim ends the proof of the push-off lemma. We prove the claim by contradiction.

{\bf Step 10.} Let us begin by indicating that arcs like $\alpha _j$ are easy to recognize. Up to isotopy fixing the boundary of $\alpha _j$ and respecting $\Gamma$, they are all like the arc $\alpha$ represented in Figure 13.

\begin{figure}[ht]    
   \centering
    \includegraphics[scale=0.4]{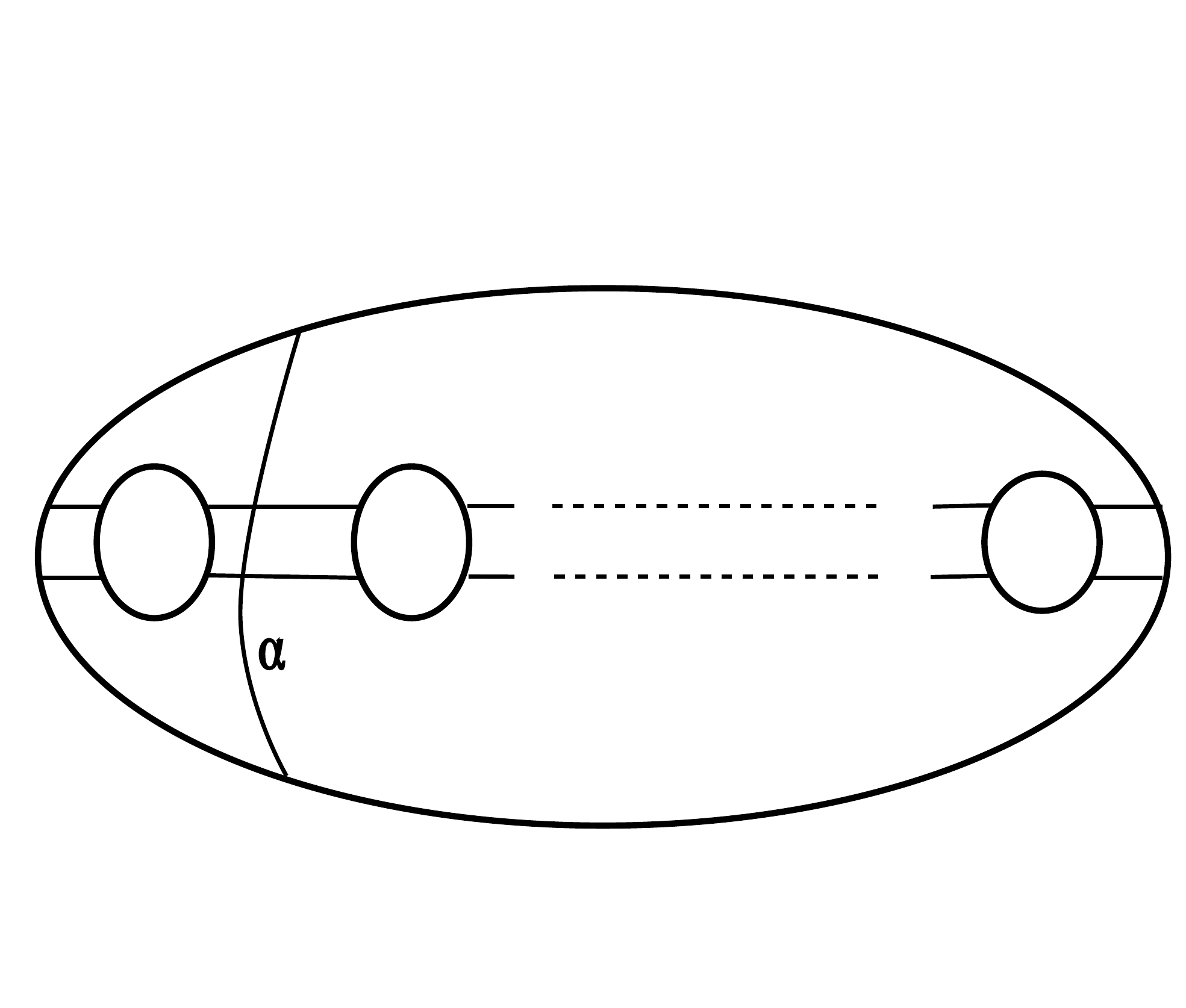}
\caption{}
\end{figure}

In Figure 13 it is easy to see that $\Delta _1$ and $\Delta _2$ cannot have the same boundary component, since $\gamma$ would cut $\Gamma$ in at least 6 points. 

{\bf Step 11.} Let us return to the beginning of Step 8. Consider $\Delta = \Delta _j$ and again let $C$ be the boundary of $\Delta$. Let $D'$ be as defined at the beginning of Step 8. Let $\{P_1, P_2\} = \alpha \cap C$. Let $D''$ be the subdiagram distinct  from $D'$ which  also has $C$ as a boundary component. The intersection $\gamma \cap D''$ consists of two arcs $\beta _1$ and $\beta _2$, one extremity of $\beta _k$ being $P_k$. Let $Q_k$ be the other extremity of $\beta _k$ ($k = 1, 2$). Both points $Q_1$ and $Q_2$ belong to the same boundary component $C'$ of $D''$ because there are only two innermost discs in $\hat {E}$. Moreover, we have $\beta _k \cap \Gamma = \emptyset$ because $\Delta _1$ and $\Delta _2$ have no boundary component in common. Since $D''$ has at least two boundary components, it cannot be a singleton.

{\bf Step 12.} Suppose that $D''$ is a band diagram. Then the existence of the  arcs $\beta _1$ and $\beta _2$ with the properties stated in Step 11 implies that the band in $D''$ extends the band in $D'$. As a consequence $C$ could be removed. This contradicts the minimality of $\mathcal {C}_m$.

{\bf Remark.} It is here that the minimality of $\mathcal {C}_m$ is used in the proof of the push-off lemma.

{\bf Step 13.} Therefore we are left with the only possibility that $D''$ is basic. We shall see that this leads to a contradiction. Let $\Delta '$ be the disc in $S^2$ bounded by $C'$ such that $D'' \cap \Delta ' = C'$. We have card$(\gamma \cap \Delta \cap \Gamma) = 2$. The circle $C'$ is broken up into two arcs $C'_1$ and $C'_2$ such that $C' = C'_1 \cup C'_2$ and $C'_1 \cap C'_2 = \{Q_1, Q_2\}$. We cannot have $C'_j \cap \Gamma =  \emptyset$ otherwise $C'_j$ would be a compressing arc for $\gamma$ (compare with step 6). Hence card$(C'_j \cap \Gamma) = 2$ for $j = 1, 2$. With $C_1$ or  $C_2$, then with $\beta _1$ and $\beta _2$ and finally with $C'_1$ or  $C'_2$ we can construct several Haseman circles in $D''$ which are not boundary parallel. This contradicts the fact that $D''$ is basic.

{\bf This concludes the proof of the claim and  the proof of the push-off lemma.}

\subsection{Uniqueness}

To complete the proof of uniqueness for minimal Conway families, we argue as follows. Let $\mathcal {C}'_m$ and $\mathcal {C}''_m$ be two of these. For each crossing point of $\Gamma$ both $\mathcal {C}'_m$ and $\mathcal {C}''_m$ contain a circle which surrounds it. It is easy to make these circles coincide by an isotopy. Therefore, from now on, we shall only consider circles in $\mathcal {C}'_m$ and $\mathcal {C}''_m$ which do not bound singletons. 
\\
Let $C$ be a circle of $\mathcal {C}''_m$. By the push-off lemma, we can isotope $C$ off $\mathcal {C}'_m$. Consider the diagrams of $\Gamma$ determined by $\mathcal {C}'_m$. Let $D'$ be the one which contains $C$. There are three possibilities.
\\
1) If $D'$ is basic, then $C$ is boundary parallel and we can make it coincide with a member of $\mathcal {C}'_m$ via an isotopy. 
\\
2) If $D'$ is a band diagram and  $C$ is boundary parallel, we again isotope it to a boundary component of $D'$, i.e. to a member of $\mathcal {C}'_m$.
\\
3) If $D'$ is a band diagram but  $C$ is not boundary parallel, we do nothing.

We then consider another member of $\mathcal {C}''_m$ and proceed as above. And so on. Eventually, we isotope each circle of $\mathcal {C}''_m$ to either:
\\
i) a circle of $\mathcal {C}'_m$, or
\\
ii) a circle in a band diagram which is not parallel to a boundary component of that diagram.

Up to now, we have only used the fact that $\mathcal {C}''_m$ is a family of Haseman circles. Since $\mathcal {C}''_m$ is a Conway family, the isotoped $\mathcal {C}''_m$ contains $\mathcal {C}'_m$. Since $\mathcal {C}''_m$ is minimal, there is no circle of the isotoped  $\mathcal {C}''_m$ satisfying ii).

{\bf This concludes the proof of the Main Theorem.}

\subsection{Final remarks} 

1) In a band diagram with at least four boundary components ($v \geq 4$) there are several Haseman circles which are not boundary parallel. For each of them, there are other Haseman circles which cut it in an essential way: their intersection cannot be removed by an isotopy. On the contrary, circles from the minimal Conway family can be characterized by the fact that their intersection with any Haseman circle is inessential: it can be removed by an isotopy. This contrast is an easy 2-dimensional reading of a much deeper result in dimension  3 due to Walter Neumann and Gadde Swarup \cite{NeSw} where incompressible tori play the role of Haseman circles.

2) The proof of the Main Theorem illustrates why maximal families of Haseman circles are not unique and why minimal Conway families are unique. It also shows that the lack of uniqueness is under control. The moral is this: there are two kinds of Haseman circles for a link projection in $S^2$. Either a Haseman circle  belongs to $\mathcal {C}_m$ or it can be isotoped into a band diagram. Haseman circles in a band diagram are easy to classify.

\section {The position of flypes}

\begin{theorem}
(Position  of flypes.) Let $\Pi$ be a link projection in $S^2$ and suppose that a flype occurs in $\Pi$. Then, its active crossing point belongs to a twisted band subdiagram determined by $\mathcal{C}_{can}$. The flype moves the active crossing point
\\
1) either to the twist to which it belongs,
\\
2) or  to another twist of the same band diagram.
\end{theorem}

Sometimes, a flype of type 1) is called an {\bf inefficient flype} while one of type 2) is called {\bf efficient}. We are interested  mainly in efficient flypes.

{\bf Comments.} If a projection is equal to its polyhedral part, then it has no flypes. It is not necessary to suppose that $\Pi$ is alternating. See Section 7 for more information about alternating diagrams.

{\bf Definition.} We call the set of crossing points of the twists of a given twisted band diagram a {\bf flype orbit}.
\\
Roughly speaking, the theorem asserts that a flype moves an active crossing point inside the flype orbit to which it belongs. Note that two distinct flype orbits are disjoint. This can be interpreted as a loose kind of commutativity of flypes. Compare with \cite{Calvo}.

{\bf Proof of Theorem .}
\\
It is obvious from their definition that, in order to locate the flypes, we have to locate the Haseman circles. This is exactly what the canonical decomposition of a link projection  into its arborescent and polyhedral parts does. The push-off lemma stated and proved in Section 3 implies that a Haseman circle for a link projection is isotopic (with respect to $\Pi$) to either:
\\
Type 1. The boundary of a singleton.
\\
Type 2. A canonical Haseman circle.
\\
Type 3. A circle contained in a twisted band diagram and not of type 1 or 2.

A canonical Haseman circle can be of three types:
\\
Type 2a. A circle separating two jewels.
\\
Type 2b. A circle separating two twisted band diagrams.
\\
Type 2c. A circle separating a jewel and a twisted band diagram.

Now Figures 3a and 3b show that in a flype two Haseman circles joined by bands are involved. A case by case examination of the three types listed above shows that these circles are isotopic to circles in a twisted band diagram and that at least one of them is not boundary parallel. When the flype is efficient, such a circle is isotopic to a circle represented by a dotted line in Figure 14.

\begin{figure}[ht]    
   \centering
    \includegraphics[scale=0.4]{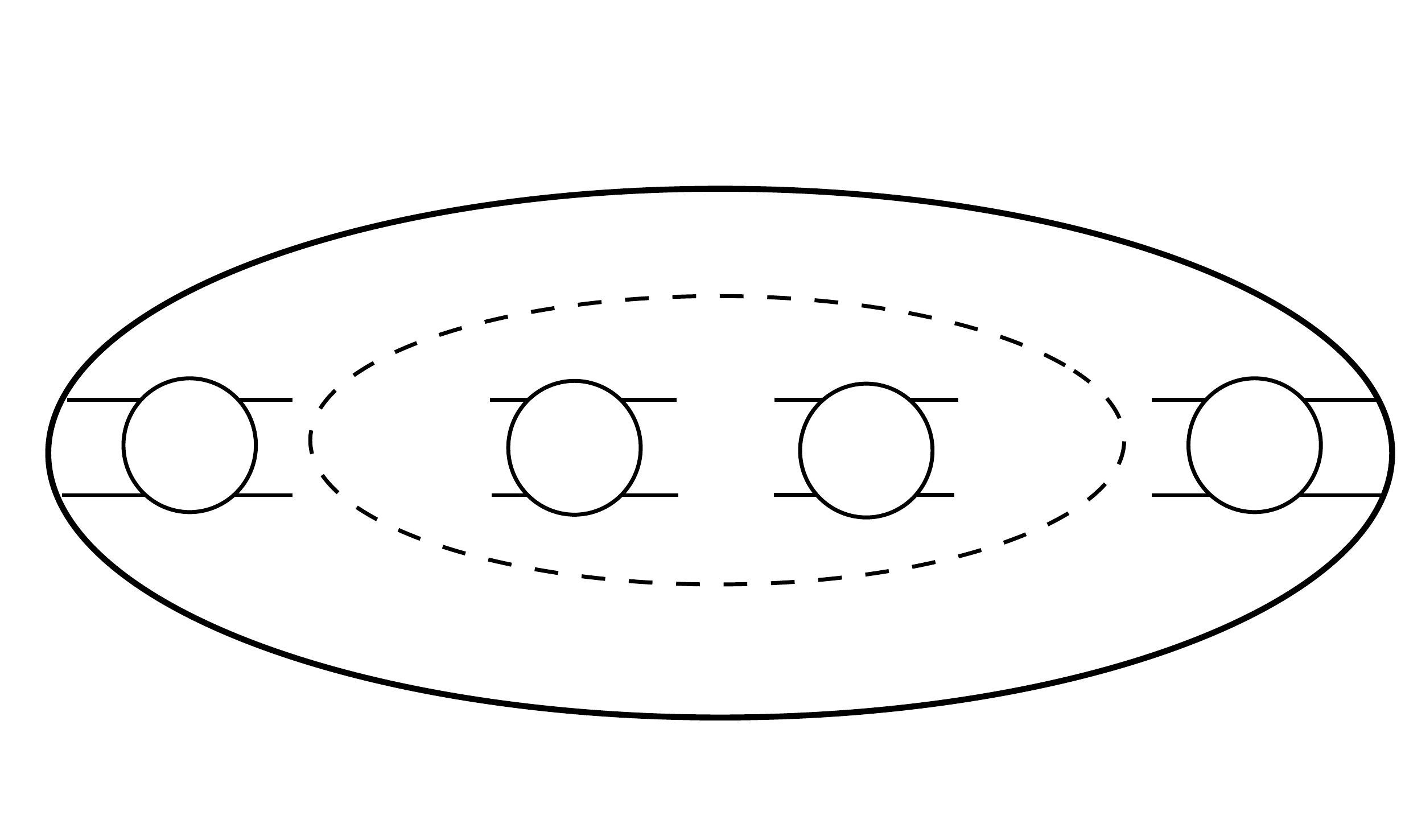}
\caption{}
\end{figure}

\section{Coding an arborescent diagram \`a la Bonahon-Siebenmann}

Let us  again consider  a twisted band diagram as pictured in Figure 15.

\begin{figure}[ht]    
   \centering
    \includegraphics[scale=0.4]{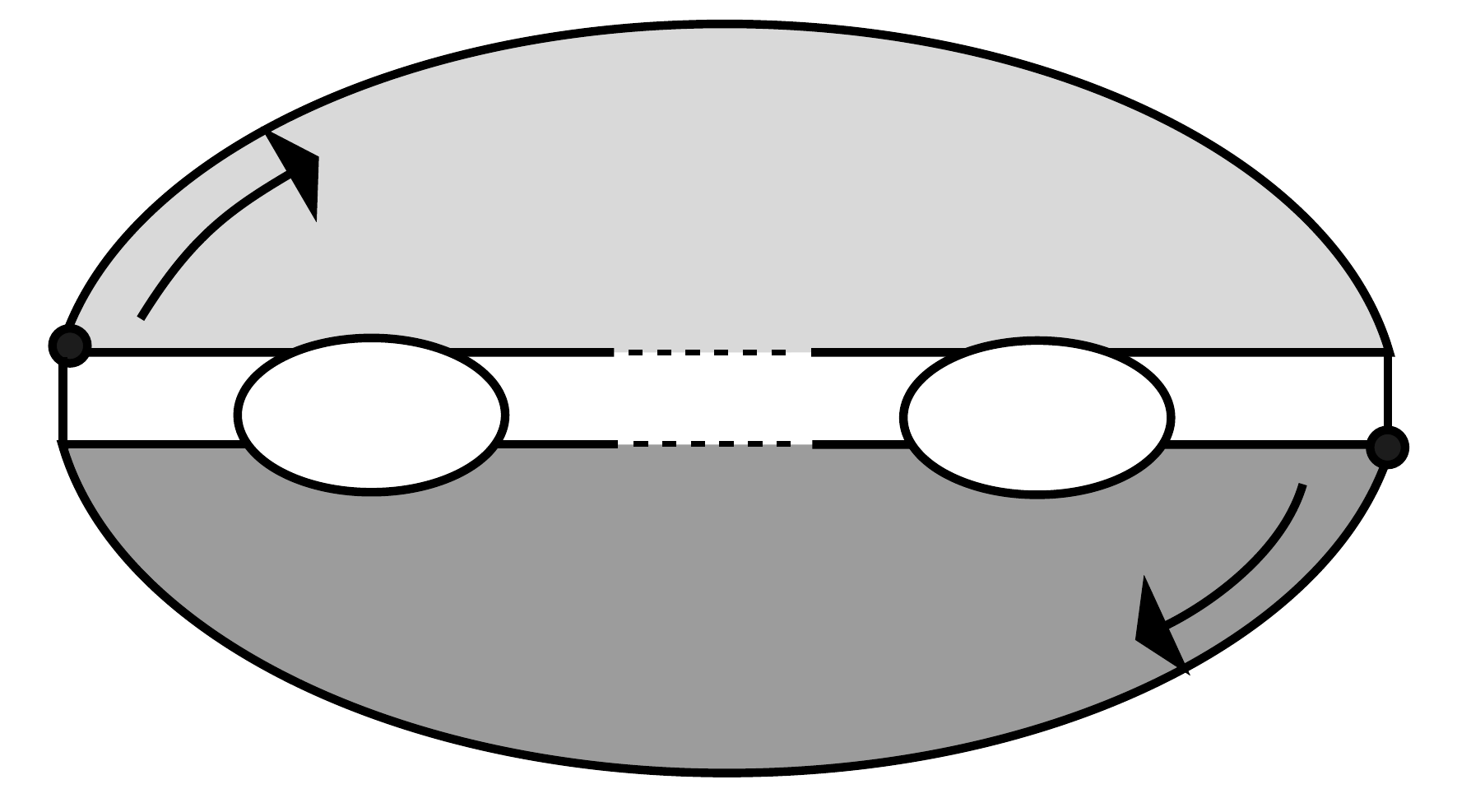}
\caption{}
\end{figure}

In this figure, two domains  have been shaded, one lightly and one darkly. To code this diagram, the idea of Bonahon-Siebenmann is:

1. To choose one of the two shaded domains.
\\
2. To travel clockwise along the boundary of the chosen domain.
\\
3. To record the events in the order  in which they are encountered during the travel. An event is either a twist or an arc of a canonical Haseman circle.

The twisted band diagram is symbolized by a planar graph constructed step by step as follows. We begin with a vertex and attach $(k+1)$ edges to it. The opposite extremity of an edge has no vertex (it is a free edge). Now, the plane is supposed to be canonically oriented. Hence, we can speak of the clockwise or counter-clockwise cyclic order around the vertex. In \cite{BoSi} the counter-clockwise order is chosen. The edges represent the canonical circles in the order in which they are met. In a sector between two succesive edges, we write the corresponding intermediate weight. Note that if we had chosen the other shaded domain in the twisted band diagram and had performed the same construction, the two graphs would differ by a reversal of the cyclic order around the vertex. We call this (elementary) graph an {\bf aster}.

Suppose then that we have an arborescent diagram $D$. The (brillant!) idea of Bonahon-Siebenmann is to code it by a planar tree (whence the name) making use of the canonical decomposition of $D$ in the subdiagrams determined by $\mathcal{C}_{can}$. A slight difficulty arises in the choice of the shaded domain for each twisted band subdiagram. We wish to make a choice for one of them, and then to propagate it. Here is how we proceed. The initial choice is represented in some twisted band diagram by a {\bf fat point} and by an arrow emanating from it, drawn a little bit inside the domain. The arrow represents the clockwise direction in which the boundary of the chosen domain is travelled. In Figure 15, the two possible starting points are shown, one for each domain. Once the starting point has been chosen, we travel along the boundary of the domain, following the arrow. In our travel, we  move alternately along an arc of a canonical circle, then along a twist, then along an arc, and so on. Each time we hit an arc of a canonical circle, we draw a fat point. It will work as starting point for the domain which lies on the other side of the canonical circle. In this manner, we propagate fat points everywhere in $D$. A conflict cannot arise, since the situation is essentially tree-like, reflecting the fact that planar surfaces have genus zero.

To complete the construction of the planar tree, we proceed as follows. Each twisted band diagram is already symbolized by an aster. If two of these meet in a common boundary component, we glue  the corresponding  asters together along the edge which represents the common boundary component. That edge is now no longer free. There is essentially one way to do the glueing in the oriented plane. At the end, we get a planar tree. The graph cannot have circuits since we work with diagrams in planar surfaces.
\begin{figure}[ht]    
   \centering
    \includegraphics[scale=0.55]{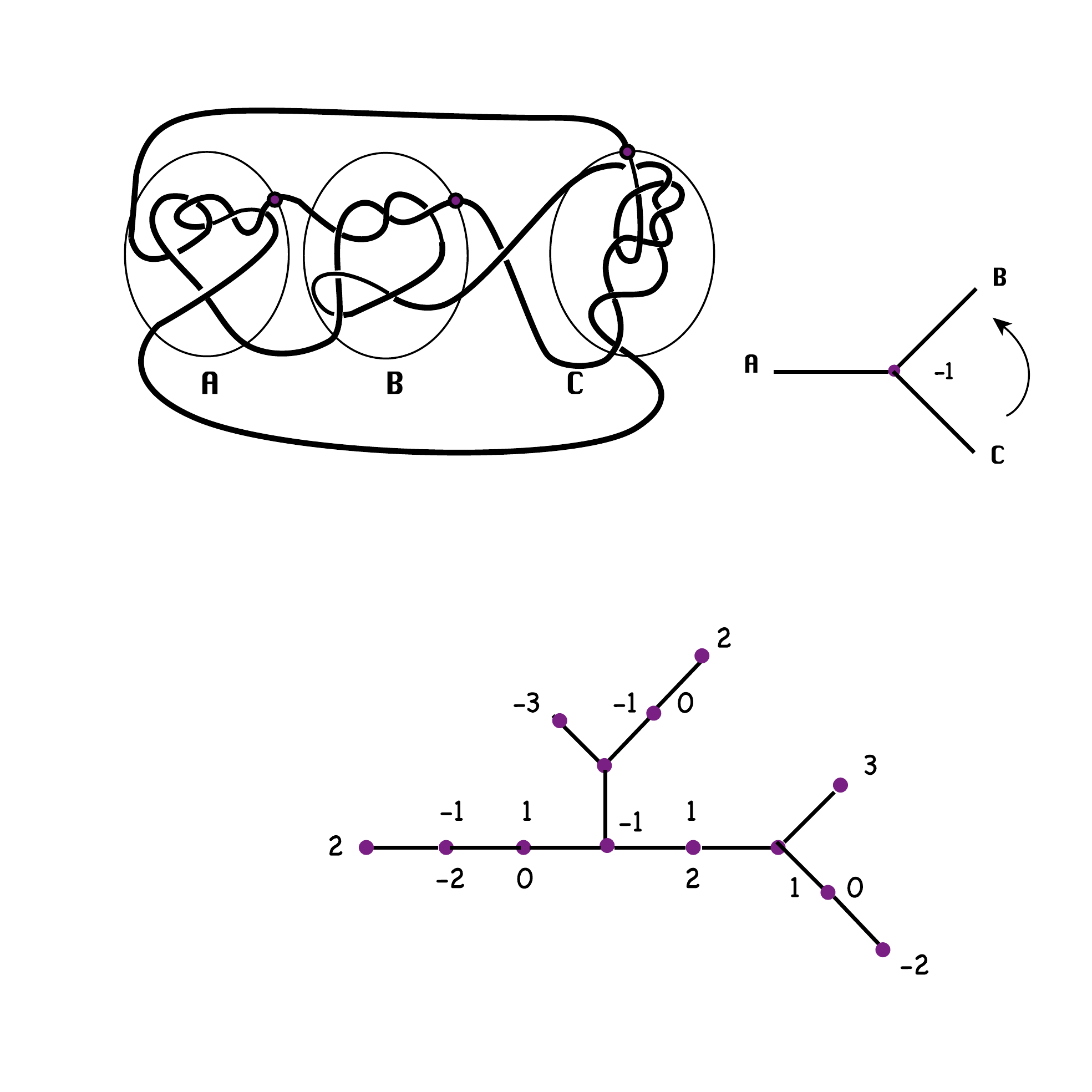}
\caption{ An arborescent knot with its weighted planar tree}
\end{figure}

{\bf Definition.} We call the tree thus obtained the {\bf planar weighted tree} associated to the arborescent diagram.

We denote it by  $\mathcal{T}$ or $\mathcal{T}(D)$. The free edges which remain represent the boundary components of $\Sigma$. In Figure 16, an arborescent diagram in $S^2$ and its corresponding planar tree are pictured; together with the canonical Haseman circles, a choice of fat points and some arrows.

{\bf Comment.}  We have called $\mathcal{T}(D)$ ``the" planar tree associated to $D$. However, the tree depends on the choice of an initial fat point. If the opposite choice is made, the tree we obtain differs from the first one by a reversal of the cyclic order around each vertex. In other words, the planar tree $\mathcal{T}(D)$ is defined only up to a reflection in a straight line. 

\section{The effect of a flype on a Bonahon-Siebenmann tree}

Let $\mathcal{T}(D)$ be the tree which codes an arborescent diagram $D = (\Sigma , \Gamma)$. Suppose that a flype occurs in $D$. We have seen in Section 4 that the active crossing point $P$ belongs to some twisted band subdiagram $D'$ of $D$. Let $V$ be the vertex of $\mathcal{T}(D)$ which represents $P$. Let $\epsilon = \pm 1$ be the sign of $P$. Let $a_i$ be the intermediate weight of the twist to which $P$ belongs. And let $a_j$ be the intermediate weight of the twist to which $P$ will move  after the flype. Because we assume the flype to be efficient we have $i \neq j$. The first effect of the flype is that $a_i$ is changed to $a'_i = a_i - \epsilon$ and $a_j$ to $a'_j = a_j + \epsilon$. Hence, the total weight is unaffected.

There is a second effect of the flype on $\mathcal{T}(D)$, related to the choice of the moving subdiagram. In the planar tree $\mathcal{T}(D)$ the intermediate weights $a_i$ and $a_j$ are written in two different sectors $A_i$ and $A_j$ around the vertex $V$. The choice of the moving subdiagram is reflected in the choice of one of the two angles with vertex $V$ which have their extremities in $A_i$ and $A_j$: we can move cyclically from $A_i$ to $A_j$ either clockwise or counter-clockwise. Let $E_1 , ... , E_s$ be the edges of the aster with vertex $V$  which belong to the chosen angle. Remove  $V$ from the tree $\mathcal{T}(D)$ and let $\mathcal{T}_i$ (for $i = 1 , ... , s)$ be the subtree which contains $E_i$. The flype affects $\mathcal{T}_i$ in the following way. Consider a vertex $W$ of $\mathcal{T}_i$ and let $d$ be the distance from $W$ to the free extremity of $E_i$ where $V$ was. Then:
\\
1) if $d$ is odd the cyclic order around $W$ is reversed,
\\
2) if $d$ is even, the cyclic order around $W$ is preserved.

This is rule F3 of Chapter 13 of \cite{BoSi}. We apologize for taking so long to present it! The proof  is by induction on $d$ and by drawing a few pictures built from models of twisted band diagrams.

{\bf Remark.} Originally, Tait used the word ``flype" to describe another move on a link projection, different from the one which is today called flype. This was  an inversion with respect to a circle in $S^2$, followed by the mirror image through $S^2$ (of course $S^2$ is the sphere into which $\Pi$ is embedded). We shall call such a move an {\bf old flype}. 

To record the action of an old flype on the tree $\mathcal{T}(D)$ there is a slight ambiguity due to the choice of the fat points (see Section 5 above). Therefore, there are two ways to draw the planar tree for the image of an arborescent diagram after an old flype. Of course, they are related by a reversal of the cyclic order at each vertex. One of these two trees is obtained from the original one as follows:
\\
1. Choose a vertex $V$ of $\mathcal{T}(D)$.
\\
2. Reverse the cyclic order around that vertex.
\\
3. Reverse the cyclic order around each vertex which is at even distance from $V$.

Again, the proof is easy. The vertex $V$ represents a twisted band subdiagram $D'$ of $D$. We have an original choice of  shaded domain for $D'$. For the image $fD$ of $D$ after the old flype  $f$ we look at the image $fD'$ of $D'$ and take as shaded domain for $fD'$ the image of the chosen shaded domain for $D'$. We propagate fat points from it. Drawing a few pictures helps.

{\bf Remark.} Let us again consider  a flype acting on a diagram $D$. We have seen in Figures 3a and 3b that there are two possible choices for the moving tangle. It is easy to see that the two choices are related by an old flype. See Figure 17. This agrees with the action of flypes (old and new) on planar trees.
\begin{figure}[ht]    
   \centering
    \includegraphics[scale=0.3]{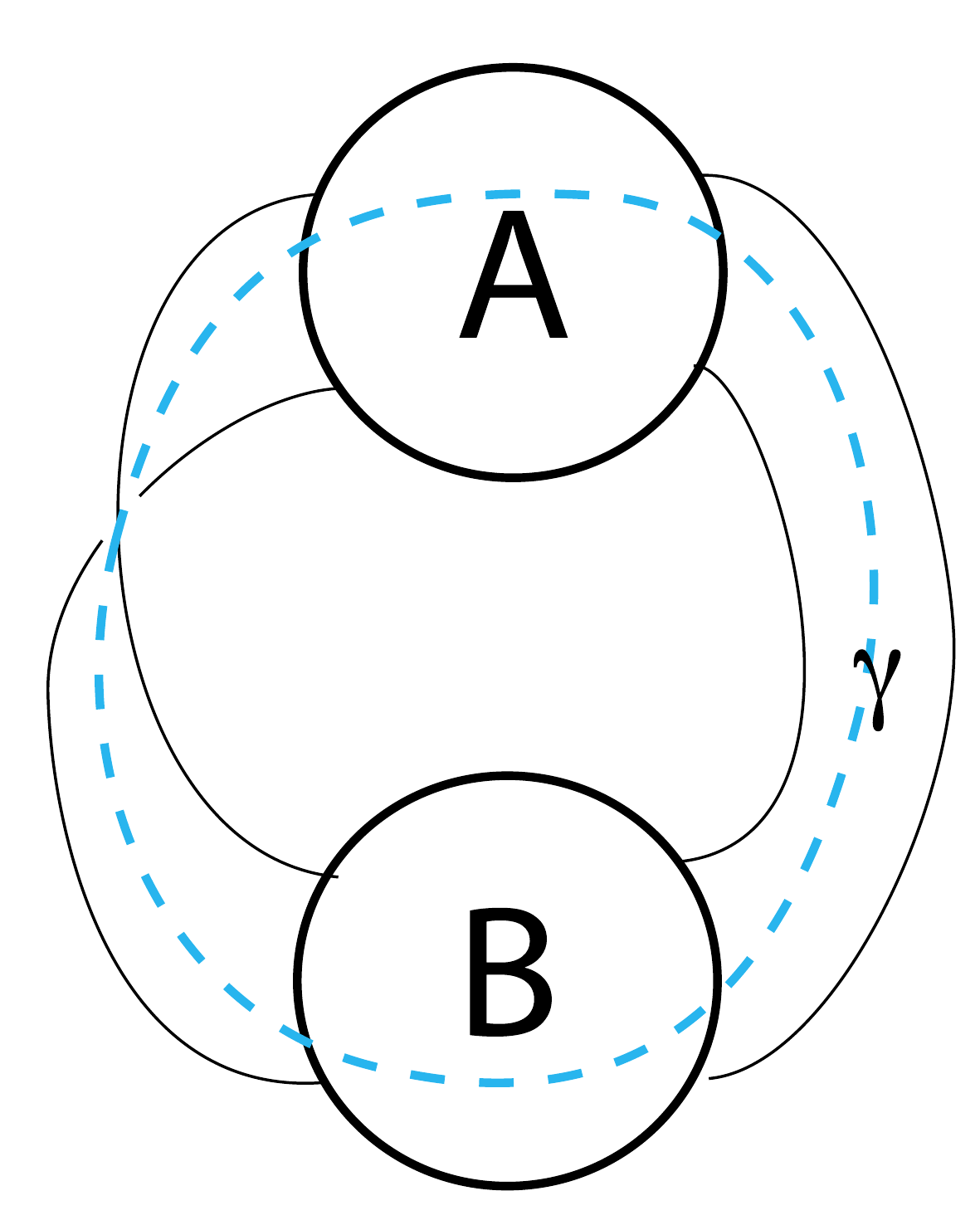}
\caption{An old flype}
\end{figure}

\section{Alternating projections}

\subsection{Arborescent alternating diagrams}

Let $D$ be an arborescent diagram and let $\mathcal{T}(D)$ be the planar tree associated to $D$. We carry out the following operations on $\mathcal{T}(D)$:
\\
1) We forget the cyclic order around each vertex.
\\
2) To each vertex we attach the total weight $a = \sum a_i$.

We call the tree obtained this way from $\mathcal{T}(D)$ the {\bf abstract tree} associated to $D$ and we denote it by $\mathcal{A}(D)$. We say that $\mathcal{T}(D)$ is above $\mathcal{A}(D)$.

{\bf Question.} When is an arborescent diagram alternating?

To answer, we  need only consider the abstract tree  $\mathcal{A}(D)$.

\begin{proposition}
An arborescent diagram $D$ is alternating iff it is possible to prescribe a sign at the vertices of $\mathcal{A}(D)$ of total weight zero in such a way that any two adjacent vertices of $\mathcal{A}(D)$ have opposite signs. 
\end{proposition}

An immediate consequence of the proposition is that the minimal crossing number of the alternating link type represented by a diagram satisfying the condition of the proposition is given by the sum of the absolute values of the weights attached to the vertices of $\mathcal{A}(D)$. See \cite{Kauffman} and \cite{Murasugi} for the proof of the Tait conjecture about the minimal crossing number of alternating links.

{\bf Proof of  Proposition 1}

1st step. Look at Figure 11 and suppose that the twisted band diagram is a subdiagram of some arborescent diagram $D$. Consider the checkerboard associated to it, together with the colour (black or white, say) for each region. Note that the two big shaded regions have the same colour and that the bands (twisted or not) have the other colour.

2nd step. Consider an edge of the planar tree $\mathcal{T}(D)$. Moving along that edge corresponds to crossing a well defined canonical Haseman circle $\gamma$. The first step implies that the bands  attached on one side of $\gamma$ have a colour opposite to that of the bands attached on the other side of $\gamma$.

3rd step. Consider a crossing point and the four sectors which are in its neighbourhood. J. B. Listing in \cite{Listing} labelled them $\lambda$ and $\delta$ according to the rule of Figure 18.

\begin{figure}[ht]    
   \centering
    \includegraphics[scale=0.3]{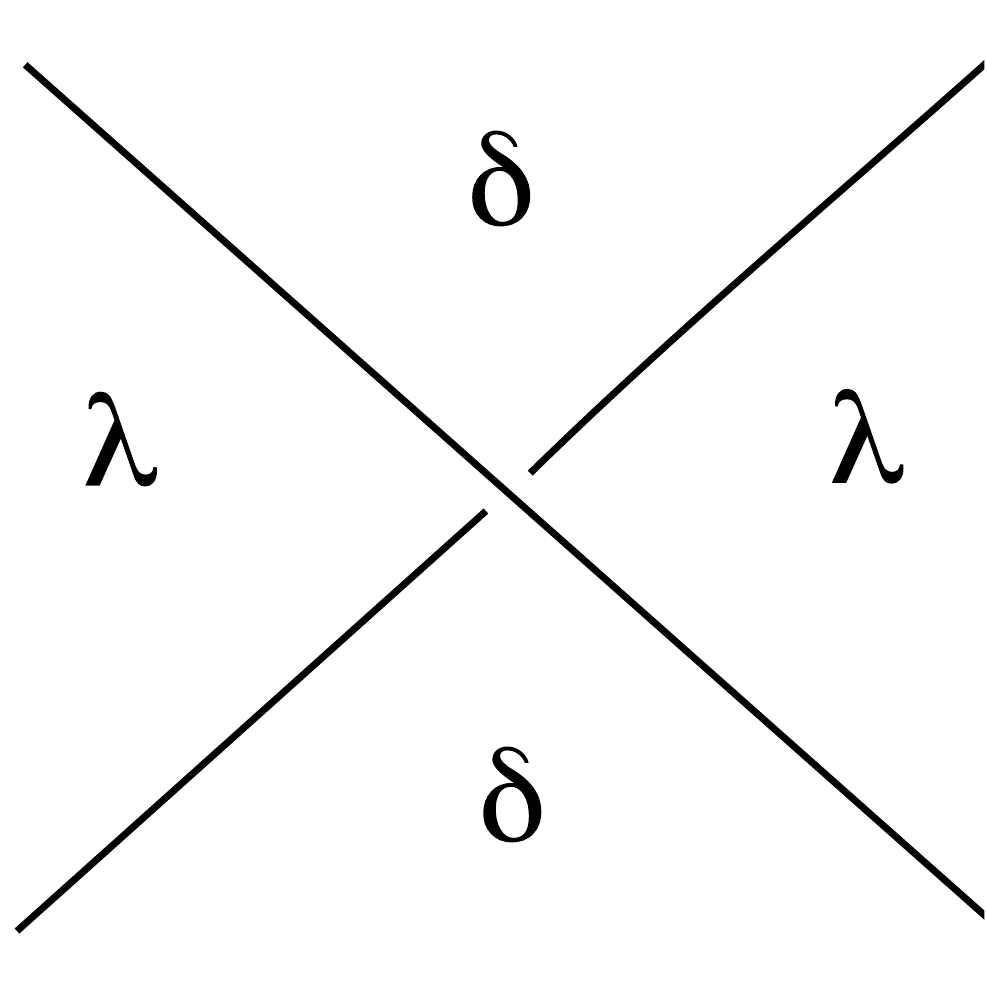}
\caption{}
\end{figure}

Listing called a region determined by a diagram a  $\bf monotyp$ region when all its sectors have the same label.

4th step.  Then consider Figure 10 where a sign is assigned to a crossing point  sitting on a band. We observe that there is the following correlation between signs and labels. If the two sectors   contained in the band have label $\lambda$  the sign is $+1$, and if their label is $\delta$  the sign is $-1$.

5th step. We remark that because of Rules 3 and 4, the sign of a crossing point in a twisted band diagram is the one given by the total weight. Therefore, we  need only  consider the abstract tree $\mathcal{A}(D)$ and not the planar tree $\mathcal{T}(D)$.

6th step. Suppose that $D$ is alternating.  Tait already observed that when a diagram is alternating then: 1) all regions determined by $D$ are monotyp, 2) the labels of the regions correspond to colours of the checkerboard. Tait's observation and the previous steps imply that the tree $\mathcal{A}(D)$ satisfies the alternating condition expressed in the proposition, because in an arborescent diagram all crossing points appear in the bands.

7th step. Conversely, suppose that $\mathcal{A}(D)$ satisfies the alternating condition for the signs. Looking at Figure 11, we quickly realize that this implies that all regions are monotyp. The 4th step implies that two adjacent regions have opposite labels. This implies that the diagram is alternating.

{\bf End of proof of Proposition 1}

{\bf Remark.} The diagram represented in Figure 16  is alternating and we see that its corresponding tree satisfies the condition of the proposition.

It is now time to recall the all-important Flyping Theorem of William Menasco and Morwen Thistlethwaite \cite{MeTh}.

{\bf The Flyping Theorem.} Let $K$  be an oriented prime alternating link type in $S^3$. Then any two minimal projections representing $K$ differ by a finite sequence of flypes (and diffeomorphisms of $S^2$).

{\bf Comments.} 
\noindent

1) The Flyping Theorem and Section 4 above imply that all minimal diagrams of an alternating arborescent link type have the same abstract tree, which of course satisfies the condition of Proposition 1.

2) Fix an abstract tree $\mathcal{A}$ satisfying the condition of Proposition 1. Consider the set of (planar equivalence classes of) planar trees which are above $\mathcal{A}$. This set is partitioned into flype equivalence classes. Two distinct classes correspond to distinct link types by the Flyping Theorem. It is interesting to study the relation which exists between two such classes. They are mutant of some sort.

\subsection{Rational diagrams}

We close this section by looking at rational diagrams, which constitute an interesting special case of arborescent alternating diagrams. Roughly speaking, they are characterized by the fact that the canonical Haseman circles are nested. A precise definition is as follows. 

{\bf Definition.} Let $D = (\Sigma , \Gamma)$ be a diagram such that:
1) $\Sigma$ is an annulus,
2) The canonical Haseman circles are boundary parallel, without taking the graph $\Gamma$ into account.
Then $D$ is a {\bf rational diagram} with two boundary components ($v = 2$).

Between two consecutive canonical circles (or between one component of $b\Sigma$ and the closest canonical circle) there are two twists. At least one of the corresponding intermediate weights is different from zero. If both are non-zero, they have the same sign, by Rule 4.

A typical picture of such a diagram with its corresponding planar tree is given in Figure 19. We also show in this figure one of the  two possible  choices of fat points and arrows. This kind of figure was often called a ``cardan"  (a ``universal joint" in English) by Larry Siebenmann and his students at Orsay around 1980. This planar tree  is reminiscent of a {\bf bamboo} with two free extremities.
\begin{figure}[ht]    
   \centering
    \includegraphics[scale=0.6]{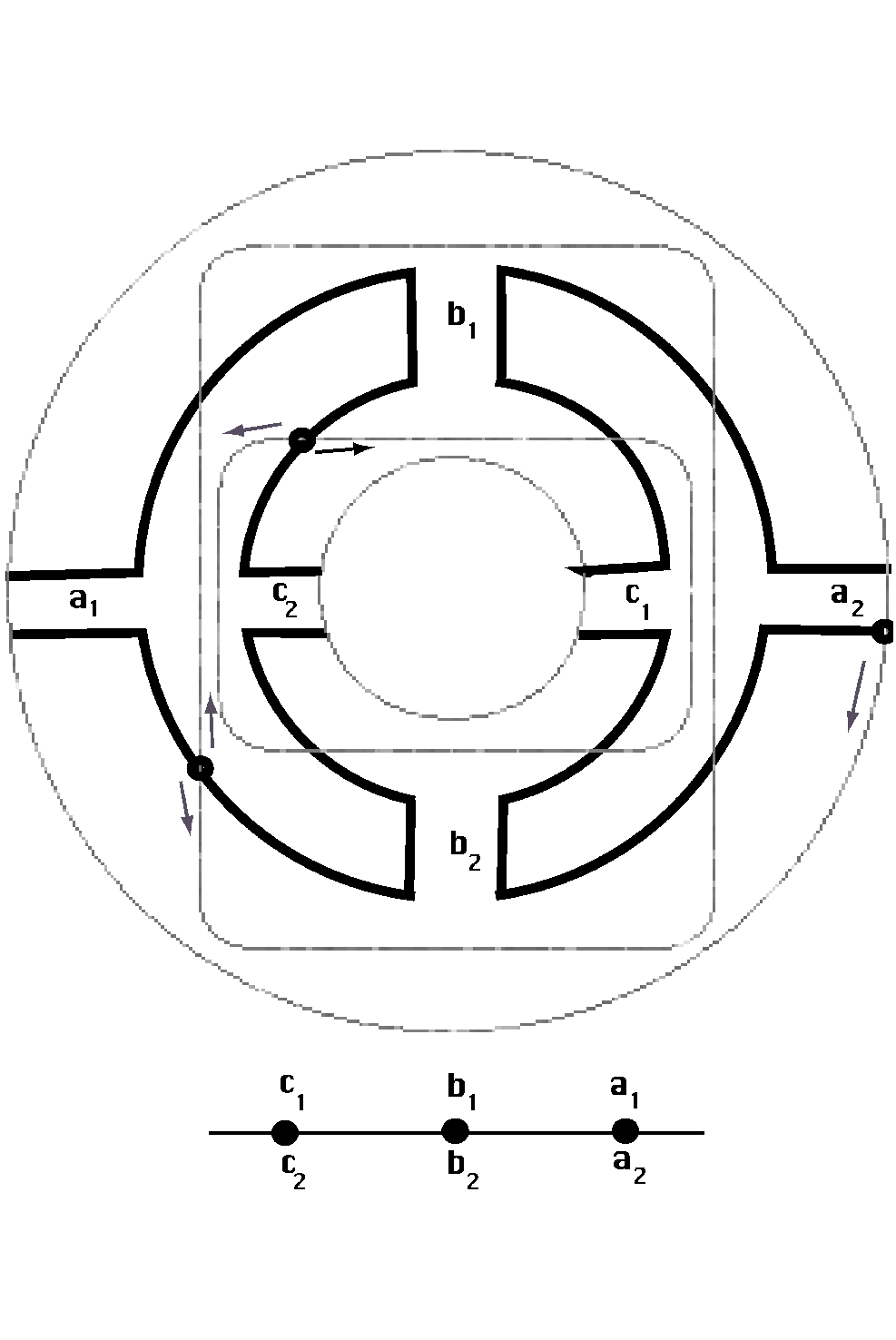}
\caption{A rational diagram with its weighted planar tree}
\end{figure}

{\bf Remark.} The building blocks of a planar tree $\mathcal{T}(D)$ are of two types: asters (with at least three edges) and bamboos. These asters correspond to twisted band diagrams (with at least three boundary components) and bamboos correspond to rational diagrams.

In a rational diagram with $v = 2$ each boundary component of  the annulus $\Sigma$ bounds a disc in $S^2$ far from $\Sigma$. If we wish, we can fill one of them, or both, with a twist in such a way that the former boundary component becomes a canonical circle. This only requires that the absolute value of the new weight be at least equal to two. We thus obtain a rational diagram with one or no boundary component. If the new $\Sigma$ is equal to $S^2$, the diagram represents a two-bridge knot or link.

Rational diagrams give rise to many flypes. The typical effect of one of them is to modify the decomposition of a given total weight into the sum of two intermediate ones. Typically $a = a_1 + a_2$ is changed into $a = (a_1 - \epsilon) + (a_2 + \epsilon)$ with $ \epsilon = \pm 1$. In the bamboo-shaped tree, the modification takes place at a vertex $V$ of valency  2. The choice of one of the possible two moving tangles is represented by the choice of one of the two sub-bamboos  which cover the bamboo and intersect exactly at $V$.

\section{Appendix: Hints for the determination of $\mathcal{C}_{can}$ and $\mathcal{T}(D)$}

{\bf Question: }How do we construct the minimal Conway family $\mathcal {C}_m$ for $\Pi$ or, better, the family $\mathcal{C}_{can}$ ? The answer lies in the first part of the proof of the Main Theorem.

i) Surround each twist of maximal length by a Haseman circle. Call such a diagram a ``spire".
\\
ii) Add new Haseman circles to that family, as long as the circles do not intersect. But do not add circles inside spires. These would anyhow not contribute  to $\mathcal{C}_{can}$ .
\\
iii) When we cannot add another member to the family, we have obtained a maximal family of Haseman circles, up to the missing ones inside spires. 
\\
iv) Remove unnecessary circles until we obtain a minimal family. This is $\mathcal{C}_{can}$. Circles which are removed are some spire boundaries and circles like the one represented by a dotted line in Figure 14. 
\\
v) During the process, several choices have been involved. The final result does not depend on them.

An example is sketched in Figure 20. 
\begin{figure}[ht]    
   \centering
    \includegraphics[scale=0.6]{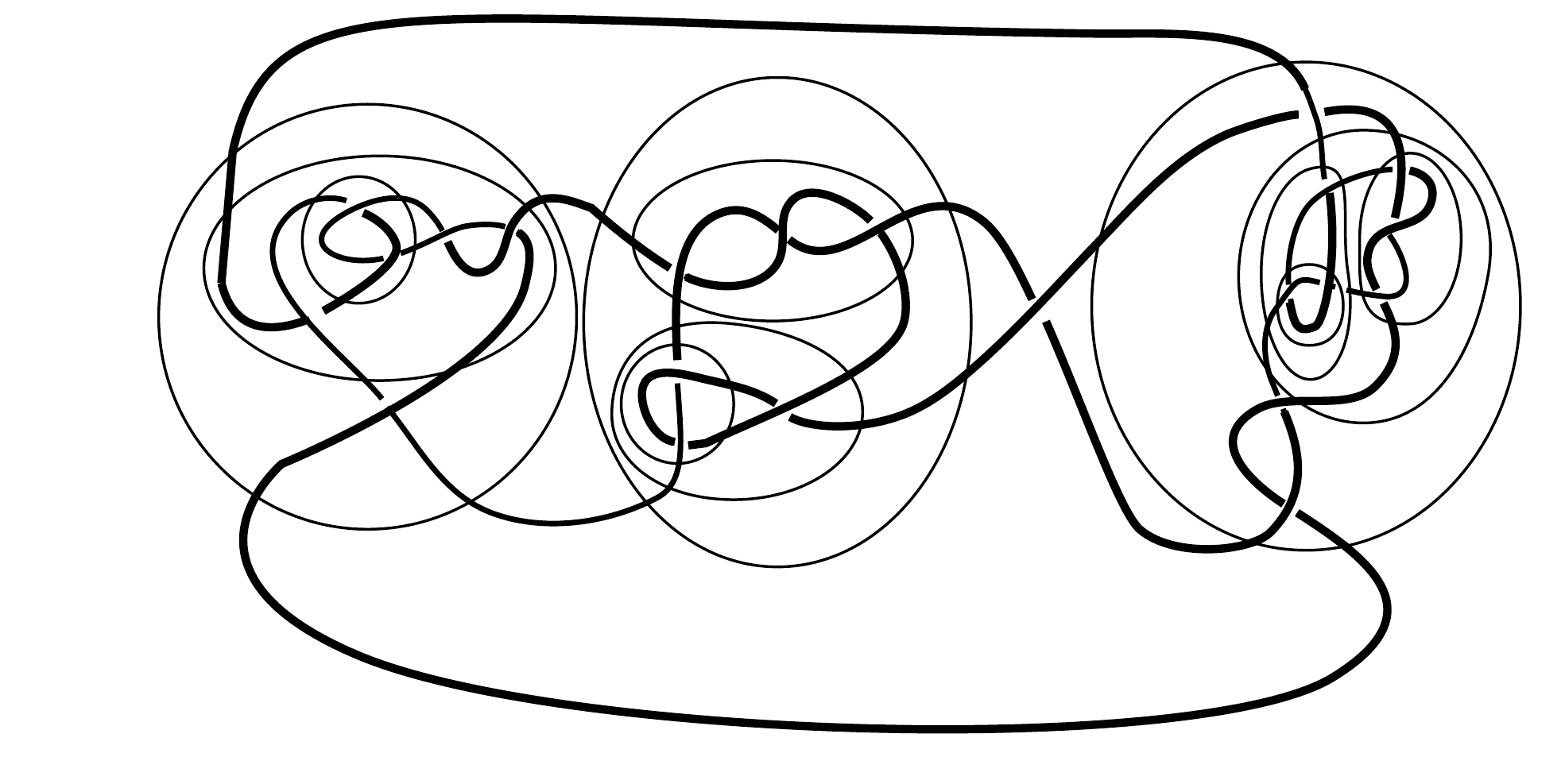}
\caption{ A knot with its minimal Conway family}
\end{figure}

Let $D$ be an arborescent diagram.
\\
{\bf Question:} How do we construct the Bonahon-Siebenmann tree $\mathcal{T}(D)$ ? 
\\
vi) Mark the location of the bands. Notice that across a circle of  $\mathcal{C}_{can}$, the position of the bands changes (this is due to  minimality!). Typically, this happens along a bamboo. In Conway's original presentation, this was specified by the position of the ``L"'s. The bands are essential to determine the signs, especially those of the singletons.
\\
vii) Proceed as explained in Section 5.

{\bf Added in proof.} The referee has rightly pointed out that Alain Caudron's thesis \cite{Ca} is related to this work. Among other things, it contains a multitude of interesting knot and link projections, some of them quite challenging. We can briefly say that Bonahon and Siebenmann \cite{Bosi} and \cite{BoSi} treat 3-dimensional objects, while we deal with 2-dimensional ones. Caudron, like Conway \cite{Conway}, works in dimension two and one half, mixing cleverly dimensions 2 and 3.

 {\bf Address}

Section de math\'ematiques
\\
Universit\'e de Gen\`eve
\\
CP  64
\\
CH-1211 GENEVE 4
\\
SWITERLAND

cam.quach@unige.ch
\\
claude.weber@unige.ch

\end{document}